\documentclass[a4paper, 11pt]{article}
\usepackage{amsmath} \usepackage{amsfonts} \usepackage{amssymb}\usepackage{latexsym}
\usepackage[english]{babel}
\usepackage{amscd}
\usepackage[dvips,final]{graphics}
\usepackage{locengli,math}
\usepackage{graphicx,pst-plot,psfig,pstricks,pst-node,pst-text,pst-tree,pst-char,pst-coil}
\newcommand{\F}{\mathcal{F}}

\newcommand{\D}{\mathbb{D}}
\newcommand{\dpe}{ \underset{\delta_2}{{\stackrel{\delta_1}{\rightleftharpoons}}}
}
\newcommand{\dpi}{
\underset{\delta_1}{{\stackrel{\tilde{\delta}_2}{\rightleftharpoons}}} }

\newcommand{\dpm}{
\underset{\delta_M}{{\stackrel{\delta}{\leftrightharpoons}}} }

\newcommand{\dpmf}{
\underset{\delta_f}{{\stackrel{\delta}{\leftrightharpoons}}} }

\newcommand{\dpmp}{
\underset{\delta}{{\stackrel{\tilde{\delta}_M}{\rightleftharpoons}}} }

\newcommand{\dpmpf}{
\underset{\delta}{{\stackrel{\tilde{\delta}_f}{\rightleftharpoons}}} }
\begin{document}
\begin{center}
\textbf{\LARGE{\textsf{From entangled codipterous coalgebras to coassociative manifolds}}}
\vskip1cm
\parbox[t]{14cm}{\large{
Philippe {\sc Leroux}}\\
\vskip4mm
{\footnotesize
\baselineskip=5mm
Institut de Recherche
Math\'ematique, Universit\'e de Rennes I and CNRS UMR 6625\\
Campus de Beaulieu, 35042 Rennes Cedex, France, pleroux@univ-rennes1.fr}}
\end{center}

\vskip1cm
{\small
\vskip1cm
\baselineskip=5mm
\noindent
{\bf Abstract:}
In our previous works, we associated with each Hopf algebra, bialgebra and 
coassociative coalgebra, a directed graph. Describing how two coassociative coalgebras, via their directed graphs, can be entangled, leads to consider special coalgebras, called codipterous coalgebras. We yield a graphical interpretation of the notion of 
codipterous coalgebra and explain the necessity to study them. By gluing these
objects, we easily obtain, 
particular objets named coassociative codialgebras and coassociative cotrialgebras and yield, thanks to their directed graph,
a simple interpretation of these structures. Similarly, we construct  Poisson algebras and dendriform coalgebras.
From these constructions, we yield a new look on coassociative coalgebras and construct an analogue of topological manifolds called coassociative manifolds. Links with non-directed graphs are also given. 
 
\tableofcontents
\section{Intoduction}
By $k$, we mean either the real field or the complex field.
Moreover, all the vector spaces will have a finite or
a denumerable basis.

In this article, we are led to manipulate several kind of coalgebras. 
These coalgebras will be visualised through a directed graph. Let us recall for the convenience of the reader, the usual definition of a directed graph.
\begin{defi}{[Directed graph]}
A {\it{directed graph}} $G$ is a quadruple \cite{Rosen}, $(G_{0},G_{1},s,t)$
where $G_{0}$ and $G_1$ are two denumerable sets respectively called the {\it{vertex set}} and the {\it{arrow set}}.
The two mappings, $s, \ t: G_1 \xrightarrow{} G_0$ are respectively called  {\it{source}} and {\it{terminus}}.
A vertex $v \in G_0$ is a {\it{source}} (resp. a {\it{sink}}) if $t^{-1}(\{v \})$ (resp. $s^{-1}(\{v\})$)
is empty. A graph $G$ is said {\it{locally finite}}, (resp. {\it{row-finite}}) if 
$t^{-1}(\{v\})$ is finite (resp. $s^{-1}(\{v\})$ is finite).
Let us fix a vertex $v \in G_0$.
Define the set $F_{v} :=\{a \in G_{1}, \ s(a)=v \}$. A {\it{weight}} associated with
the vertex $v$ is a mapping $w_v: F_{v} \xrightarrow{} k$.
A directed graph equipped with a family of weights $w := (w_v)_{v \in G_0}$ 
is called a {\it{weighted graph}}.
\end{defi}
In the sequel, directed graphs will be supposed locally finite and row finite.
Let us introduce particular coalgebras named $L$-coalgebras \footnote{This notion has been introduced in
\cite{Coa} and developed in \cite{Coa}\cite{codialg1}\cite{perorb1}\cite{tresse}.} and explain why this notion is interesting.
\begin{defi}{[$L$-coalgebra]}
A {\it{$L$-coalgebra}} $G$ 
over a field $k$ is a $k$-vector space equipped with a right coproduct,
$\Delta: G \xrightarrow{} G^{\otimes 2}$ and a left coproduct, $\tilde{\Delta}: G \xrightarrow{} G^{\otimes 2}$, verifying the
coassociativity breaking equation 
$(\tilde{\Delta} \otimes id)\Delta = (id \otimes \Delta)\tilde{\Delta}$.
If $\Delta = \tilde{\Delta}$, the coalgebra is said  {\it{degenerate}}. A $L$-coalgebra may have two counits, the right counit  
$\epsilon: G \xrightarrow{} k$, verifying $ (id \otimes \epsilon)\Delta = id$
and the left counit $\tilde{\epsilon}: G \xrightarrow{} k$, verifying $ ( \tilde{\epsilon} \otimes id)\tilde{\Delta} = id. $ A $L$-coalgebra is
said {\it{coassociative}} if its two coproducts are coassociative, in this case the equation $(\tilde{\Delta} \otimes id)\Delta = (id \otimes \Delta)\tilde{\Delta}$ is called the {\bf{entanglement equation}}.
Denote by $\tau$, the {\it{transposition}} mapping, i.e. $G^{ \otimes 2} \xrightarrow{\tau} G^{ \otimes 2}$ such that $\tau(x \otimes y) = y \otimes x$, for all $x,y \in G$.
The $L$-coalgebra $G$ is said to be {\it{$L$-cocommutative}} if for all $v \in G$, $(\Delta - \tau\tilde{\Delta})v = 0$. 
\end{defi}

Let $G$ be a directed graph equipped with a family of weights $(w_v)_{v \in G_0}$. Let us consider the free vector space $kG_0$ generated by $G_0$.
The set $G_1$ is then viewed as a subset of $(kG_0)^{\otimes 2}$ by identifying $a \in G_1$ with $s(a) \otimes t(a)$. The 
mappings source and terminus are then linear mappings still called 
source and terminus $s,t: \ (kG_0)^{\otimes2} \xrightarrow{} kG_0$, such that $s(u \otimes v) = u$ and $t(u \otimes v) = v$, for all $u,v \in G_0$. The family of weights $(w_v: F_v \xrightarrow{} k)_{v \in G_0}$ is then viewed as a family of linear mappings.
Let $v \in G_0$. Define the right coproduct $\Delta: kG_0 \xrightarrow{} (kG_0)^{\otimes 2}$, such that
$\Delta(v) := \sum_{i: a_i \in F_v} \ w_v(a_i) \ v \otimes t(a_i)$ and the left coproduct
$\tilde{\Delta}: kG_0 \xrightarrow{} (kG_0)^{\otimes 2}$, such that $\tilde{\Delta}(v) := \sum_{i: a_i \in P_v} \ w_{s(a_i)}(a_i) \ s(a_i) \otimes v$, where $P_{v} $ is the set $\{a \in G_{1}, \ t(a)=v \}$. 
With these definitions the vector space $kG_0$ is a $L$-coalgebra called a finite Markov $L$-coalgebra since its coproducts $\Delta$ and  $\tilde{\Delta}$ verify
the
coassociativity breaking equation 
$(\tilde{\Delta} \otimes id)\Delta = (id \otimes \Delta)\tilde{\Delta}$. This particular coalgebra is called in addition finite Markov ($L$-coalgebra)
because for all $v \in G_0$, the sets $F_v$ and $P_v$ are finite and the coproducts are of the form $\Delta(v) := v \otimes \cdots $ and $\tilde{\Delta}(v) := \cdots \otimes v $.

Assume that we start with the Markov $L$-coalgebra just described and
associate
with each tensor product $\lambda x \otimes y$, where $\lambda \in k$ and $x,y \in G_0$, appearing in the definition of the coproducts, a directed arrow
$x \xrightarrow{\lambda}y$. The directed graph so obtained, called the {\it{geometric support}} of this $L$-coalgebra, is up to a graph isomorphism \footnote{A {\it{graph isomorphism}} $f: G \xrightarrow{} H$ between two graphs $G$ and $H$ is a pair of bijection $f_0: G_0 \xrightarrow{} H_0$ and $f_1: G_1 \xrightarrow{} H_1$ such that
$f_0(s_G(a))=s_H(f_1(a))$ and $f_0(t_G(a))=t_H(f_1(a))$ for all $a \in G_1$. All the directed graphs in this formalism will be considered up to a graph isomorphism.}, the directed graph we start with. Therefore, general $L$-coalgebras generalise the notion of directed graph. If $G$
is a $L$-coalgebra generated as a $k$-vector space by a set $G_0$, then its geometric support 
$Gr(G)$ is a directed graph with vertex set $Gr(G)_0=G_0$ and with arrow set $Gr(G)_1$, the 
set of tensor products $a \otimes b$, with $a,b \in G_0$, appearing
in the definition of the coproducts of $G$. As a coassociative coalgebra is a particular $L$-coalgebra, we naturally construct its directed graph.
We draw attention to the fact that a directed graph can be the
geometric support of different $L$-coalgebras.
\begin{exam}{}
The directed graph:
\begin{center}
\includegraphics*[width=3.5cm]{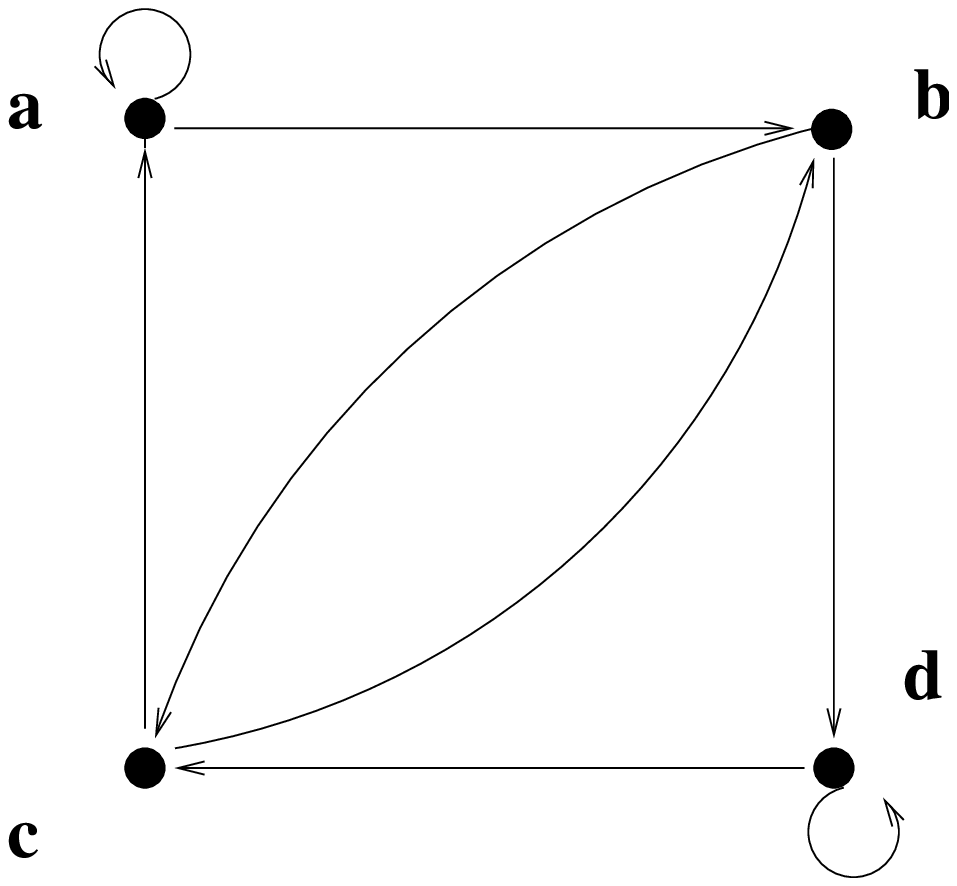}
\end{center}
\end{exam}
is the geometric support associated with the degenerate $L$-coalgebra or coassociative coalgebra, generated by $a,b,c$ and $d$, as a $k$-vector space, and
described by the following coproduct:
$
\Delta a = a \otimes a + b \otimes c, \ \
\Delta b = a \otimes b + b \otimes d, \ \
\Delta c = d \otimes c + c \otimes a, \ \
\Delta d = d \otimes d + c \otimes b
$
and the geometric support of the finite Markov $L$-coalgebra, generated by $a,b,c$ and $d$, as a $k$-vector space, and described by the right coproduct:
$
\Delta_M a = a \otimes (a + b), \ \
\Delta_M b = b \otimes (c + d), \ \
\Delta_M c = c \otimes  (a +b), \ \
\Delta_M d = d \otimes (c + d)
$
and the left coproduct: 
$
\tilde{\Delta}_M a = (a+c) \otimes a,  \ \
\tilde{\Delta}_M b = (a +c) \otimes b , \ \
\tilde{\Delta}_M c = (b + d) \otimes c , \ \
\tilde{\Delta}_M d = (b+ d) \otimes d.
$
\NB
Let $G$ be a finite Markov $L$-coalgebra. If the family of weights
used for describing right and left
coproducts take values into $\mathbb{R}_+$ and
if the right counit $\epsilon: v \mapsto 1$ exists, then the geometric
support associated with $G$
is a directed graph equipped with a family of probability vectors. 

Before going on, let us interpret what represents the $L$-cocommutativity
in the case of a finite Markov $L$-coalgebra. A directed graph is said bi-directed if for any arrow from a vertex $v_1$ to a vertex $v_2$, there exists an arrow from $v_2$ to $v_1$. To take into account the bi-orientation of a directed graph in an algebraic way, we first embed this directed graph into the finite
Markov $L$-coalgebra described above. We then notice that a directed graph is bi-directed if and only if $\Delta = \tau \tilde{\Delta}$. Therefore, in this algebraic framework, we are naturally
led to consider the $L$-cocommutator space $\ker (\Delta - \tau \tilde{\Delta})$. Dualizing this formula leads to consider an algebra $D$ equipped with two products $\vdash$ and $\dashv$ and to consider
the particular commutator space $[x,y]:= x \dashv y - y \vdash x$. The bracket, $[ -,z]$, verifies the ``Jacobi identity'', i.e. 
$[[x,y ] ,z ] = [[ x,z ] ,y ] + [ x,[ y  ,z ]],$
if $D$ is an algebra called an associative dialgebra \cite{Loday}. 

Another motivation concerning associative dialgebras is the following.
In a long-standing project whose ultimate aim is to study periodicity phenomena in algebraic $K$-theory,
J-L. Loday in \cite{Loday}, and J-L. Loday and M. Ronco in  \cite{LodayRonco} introduce several kind of algebras, one of which is 
the ``non-commutative Lie algebras'', called {\it{Leibniz algebras}}. Such algebras $D$ are described by a bracket $[ -,z]$ verifying the Leibniz identity:
$$[[x,y ] ,z ] = [[ x,z ] ,y ] + [ x,[ y  ,z ]].$$
When the bracket is skew-symmetric, the Leibniz identity becomes the
Jacobi identity and the Leibniz algebra turns out to be a Lie algebra.
A way to construct such Leibniz algebra is to start from an {\it{associative dialgebra}}, that is a $k$-vector space $D$ equipped with two associative products,
$\vdash$ and $\dashv$, such that for all $x,y,z \in D$ 
\begin{enumerate}
\item{$x \dashv (y \dashv z) = x \dashv (y \vdash z),$}
\item{$(x \vdash y) \dashv z = x \vdash (y \dashv z),$}
\item{$(x \dashv y) \vdash z = (x \vdash y) \vdash z.$}
\end{enumerate}
The associative dialgebra is then a Leibniz algebra by defining the bracket
$[x,y ] := x \dashv y - y \vdash x$, for all $x,y \in D$. 
The operad associated with associative dialgebras is then Koszul dual to the operad associated with dendriform algebras,
a {\it{dendriform algebra}} $Z$ being a $k$-vector space equipped with
two binary operations,
$ \prec \ , \ \succ: \ Z \otimes Z \xrightarrow{} Z,$
satisfying the following axioms:
\begin{enumerate}
\item {$(a \prec b) \prec c = a \prec (b \prec c) + a \prec (b \succ c),$ }
\item {$(a \succ b) \prec c = a \succ (b \prec c),$ }
\item {$(a \prec b) \succ c + (a \succ b) \succ c = a \succ (b \succ c).$ }
\end{enumerate}
This notion dichotomizes the notion of associativity since the product
$a *b = a \prec b + a \succ b$, for all $a,b \in Z$ is associative. 
Before continuing,
let us recall the following proposition from \cite{Coa}, 
\begin{prop}
Any associative dialgebra can be viewed as a dendriform algebra. 
\end{prop}
\Proof
Let $(D,\vdash, \dashv)$ be an associative dialgebra. Let $a,b \in D$.
The relations,
$ a \prec b = a \dashv b$ and $a \succ b = a \vdash b - a \dashv b,$
embed $D$ into a dendriform algebra. We notice that
$a *b := a \prec b + a \succ b = a \vdash b$ is associative. For instance,
$(a \prec b) \prec c = a \prec (b * c)$ means $(a \dashv b) \dashv c = a \dashv (b \vdash c)$ and so forth.
\eproof

\noindent
By dualizing these notions, we can easily define {\it{coassociative codialgebras}} and {\it{dendriform coalgebras}}. 
\begin{defi}{[Coassociative codialgebra]}
A {\it{coassociative codialgebra}} $D$ is a $k$-vector space equipped with
two coproducts $\delta, \hat{\delta}: \ D \xrightarrow{}D^{\otimes 2}$, verifying:
\begin{enumerate}
\item {$\delta$ and $\hat{\delta}$ are coassociative,}
\item{$(id \otimes \hat{\delta})\hat{\delta} = (id \otimes \delta) \hat{\delta} $,}
\item{$(\delta \otimes id) \delta = ( \hat{\delta} \otimes id) \delta$,}
\item{$(\delta \otimes id) \hat{\delta} = ( id \otimes \hat{\delta}) \delta$.}
\end{enumerate}
\end{defi}
\begin{defi}{[Dendriform coalgebra]}
A {\it{dendriform coalgebra}} $Z$ is a $k$-vector space equipped with
two coproducts $\delta, \hat{\delta}: \ Z \xrightarrow{} Z^{\otimes 2}$, verifying:
\begin{enumerate}
\item{$(id \otimes (\delta + \hat{\delta}))\hat{\delta} = (\hat{\delta} \otimes id) \hat{\delta} $,}
\item{$(id \otimes \hat{\delta}) \delta = ( \delta \otimes id) \hat{\delta}$,}
\item{$((\hat{\delta} + \delta) \otimes id) \delta = ( id \otimes \delta) \delta$.}
\end{enumerate}
This notion dichotomizes the notion of coassociativity since 
$(\delta + \hat{\delta})$ is a coassociative coproduct.
\end{defi}
Similarly, so as to define a {\it{``non commutative version''}} of Poisson algebra, J-L. Loday and M. Ronco introduce in \cite{LodayRonco}, the
{\it{associative trialgebras}}. Let us just mention that the operad associated with trialgebras is Koszul dual to the operad associated with dendriform trialgebras, dendriform trialgebras being $k$-vector spaces equipped with three laws $\prec$, $\succ$, $\cdot$, verifying special axioms, see also \cite{Lertribax}. Similarly, the law $*$ such that $x *y := x \prec y + x \succ y + x \cdot y$, will be associative.
Here, we are interested in contructions of
coassociative cotrialgebras. 
\begin{defi}{[Coassociative cotrialgebra]}
\label{cotri}
A {\it{coassociative cotrialgebra}} $T$ is a $k$-vector space equipped with
three coproducts $\Delta, \delta, \hat{\delta}: \ T \xrightarrow{} T^{\otimes 2}$, verifying:
\begin{enumerate}
\item {$\Delta$ is coassociative,}
\item {$(T, \delta, \hat{\delta})$ is a coassociative codialgebra,}
\item{$ (\hat{\delta} \otimes id)\hat{\delta} = (id \otimes \Delta)\hat{\delta},$}
\item{$ (\Delta \otimes id)\hat{\delta} = (id \otimes \hat{\delta})\Delta,$}
\item{$ (\hat{\delta} \otimes id)\Delta = (id \otimes \delta)\Delta,$}
\item{$ (\delta \otimes id)\Delta = (id \otimes \Delta)\delta,$}
\item{$ (\Delta \otimes id)\delta = (id \otimes \delta)\delta.$}
\end{enumerate}
\end{defi}

Let us notice that all these vector spaces are equipped with two coproducts
verifying the entanglement equation. Therefore, there exists
a graphical representation of these objets determined by
directed graphs. These directed graphs will be the
geometric supports associated with these particular $L$-coalgebras.

In addition to a graphical representation of the axioms defined above, the main results of this article lie in the constructions of
$L$-Hopf algebras, i.e. $L$-bialgebras equipped with two linear maps $\sigma$ and $\tilde{\sigma}$ such that $m(id \otimes \tilde{\sigma}) \tilde{\Delta}  := 1\tilde{\epsilon}$ and $m(\sigma \otimes id) \Delta  := 1\epsilon$.
We construct also 
coassociative codialgebras as well as  coassociative cotrialgebras. For that, we start with two coassociative coalgebras and construct bridges, i.e. coproducts, between them. 
Graphically speaking, the r\^ole of these bridges is to establish a connection between the two coalgebras.  Therefore,
we are led to consider a vector space equipped with two coproducts,
a coassociative one and a bridge, i.e. an extra left (or right) comodule on itself. Such a vector space will be called a codipterous (or anti-codipterous) coalgebra. 
\begin{center}
\includegraphics*[width=7cm]{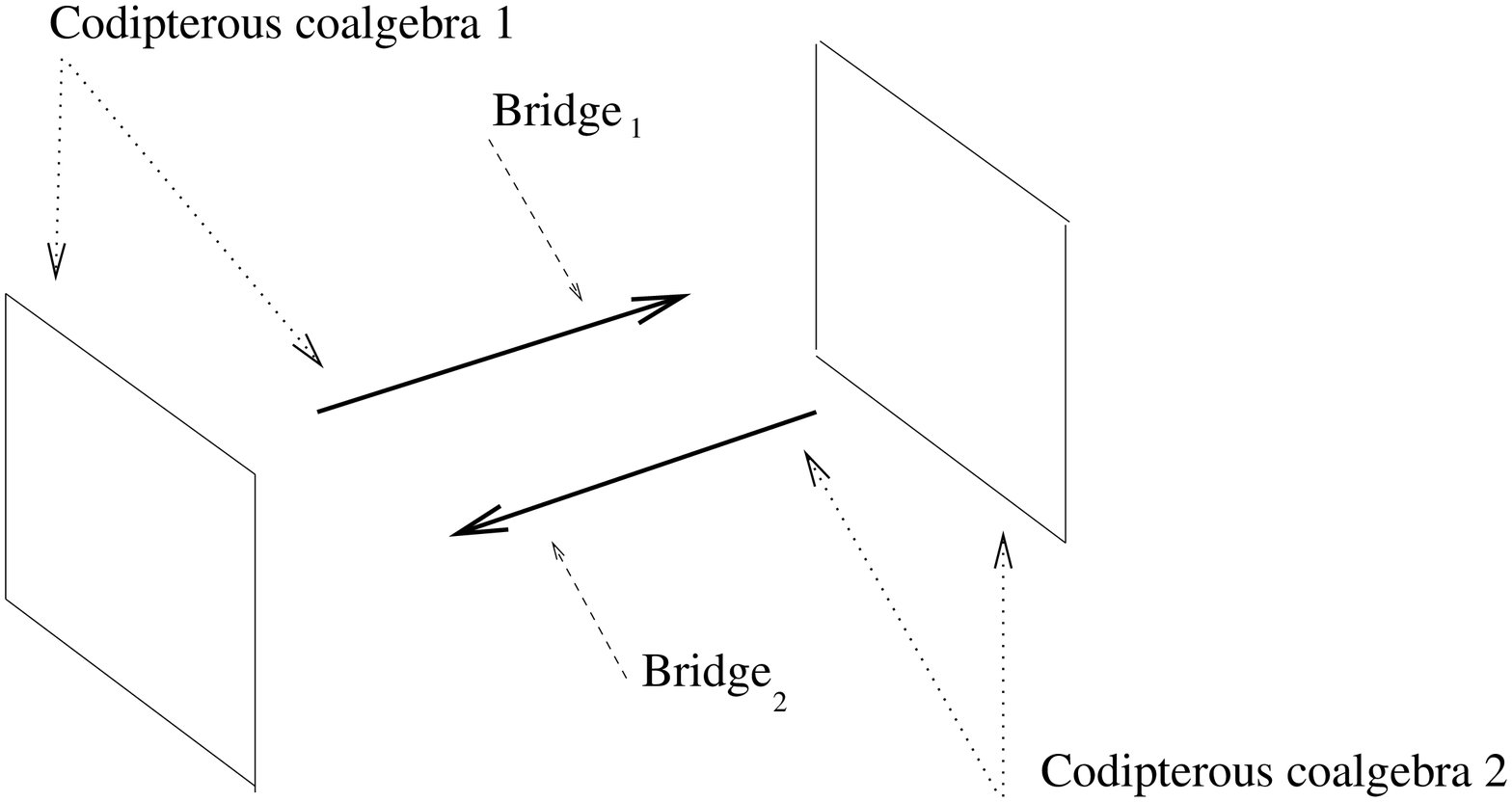}

\textbf{\begin{scriptsize} Entanglement of two codipterous coalgebras. \end{scriptsize}}
\end{center}
Schematically speaking, a codipterous coalgebra is a face represented by a square symbolizing a coassociative coalgebra and by a bridge, symbolizing the left comodule.  We realize an entanglement 
by ``gluing'' two codipterous coalgebras, such that their bridges verify the entanglement equation.
There are two important interpretations of such a construction. Firstly, suppose that we consider a $k$-vector space $E$, with dimension large enough to include several codipteous coalgebras. Suppose we view 
a codipterous coalgebra like an atom in chemistry and the $k$-vector space $E$ as a chemical solution. The entanglement of two codipterous coalgebras will yield a so-called $L$-molecule. By entangling several codipterous coalgebras, we can construct more and more complicated $L$-molecules. 
Secondly, what we have obtained in the construction of such a $L$-molecule, can be viewed as an interlocking of several coassociative coalgebras. Graphically, as we associate a directed graph with each coalgebra, 
a $L$-molecule is viewed as a special directed graph which admits a ``coassociative'' covering. By analogy with the topological manifolds, which are covered with open sets, we define coassociative manifolds. Indeed, the directed graph associated with a coassociative manifold will admit a covering with directed graphs, geometric support of coassociative coalgebras. Hence, the r\^ole of an open set will be played by a coassociative coalgebra.

Let us briefly introduce the organisation of the paper.
In section 2, we define the notion of self-entanglement and construct dialgebras, trialgebras, Poisson and Leibniz algebra and dendriform algebras. In section 3 and 4, we construct other types of entanglement which lead
to $L$-Hopf algebras and end this paper to the
notion of coassociative manifold and $L$-molecule.
We show links between non-directed graphs and Markovian coassociative manifolds.
\section{Entanglement of codipterous coalgebras}
The aim of this section is to construct associative dialgebras, associative trialgebras, Poisson and Leibniz algebras and $L$-Hopf algebras. Let us start with the definition of a codipterous coalgebra and the notion of pre-dendriform coalgebra. 
\begin{defi}{[codipterous coalgebra, anti-codipterous coalgebra]}
A $k$-vector space $\mathbb{D}$ equipped with two coproducts
$\Delta, \delta: \D \xrightarrow{} \D^{\otimes 2}$ verifying:
\begin{enumerate}
\item{Coas: $(\Delta \otimes id)\Delta = (id \otimes \Delta)\Delta$.}
\item{Codip: $(\Delta \otimes id)\delta = (id \otimes \delta)\delta$.}
\end{enumerate} 
is called a {\it{codipterous coalgebra}} \footnote{A codipterous coalgebra (resp. anti-codipterous coalgebra) can also be viewed as a coassociative coalgebra with an extra left (resp. right) comodule on itself.}. We call {\it{bridge}}, the coproduct $\delta$.
Similarly, we call an {\it{anti-codipterous coalgebra}} a $k$-vector space $\hat{\mathbb{D}}$ equipped with two coproducts
$\Delta, \hat{\delta}: \D \xrightarrow{} \D^{\otimes 2}$ such that $\Delta$
is coassociative and
$(id \otimes \Delta)\hat{\delta} = (\hat{\delta} \otimes id )\hat{\delta}$.
\end{defi}
\begin{defi}{[Entanglement of codipterous coalgebras]}
Let $\D_1$ and $\D_2$ be two sub-spaces of a vector space $E$.
Suppose $(\D_1, \Delta_1, \delta_1)$ and $(\D_2, \Delta_2, \delta_2)$ are two codipterous coalgebras. The spaces
$\D_1$ and $\D_2$ are said {\it{entangled}} if the bridges
$\delta_1, \ \delta_2$ verify the entanglement equation. The entanglement is said {\it{achiral}} if the entanglement equation is verified whatever the positions occupied by $\delta_1$ and $ \delta_2$ and {\it{chiral}} otherwise. The achiral (resp. chiral) entanglement is denoted by $[\D_1 \stackrel{\delta_1 , \ \delta_2}{\longleftrightarrow} \D_2]$ (resp. $[\D_1 \underset{\delta_2}{
{\stackrel{\delta_1}{\rightleftharpoons}}}\D_2]$, if $(id \otimes \delta_2)\delta_1 = (\delta_1 \otimes id)\delta_2$ and $[\D_1 \underset{\delta_1}{
{\stackrel{\delta_2}{\rightleftharpoons}}}\D_2]$, if $(id \otimes \delta_1)\delta_2 = (\delta_2 \otimes id)\delta_1$). 

Similarly, a codipterous coalgebra $(\D, \Delta, \delta_1)$ is entangled with an anti-codipterous coalgebra $(\D, \Delta, \hat{\delta}_1)$ if the bridges $\delta_1$ and $\hat{\delta}_1$ verify the entanglement equation. \textsf{More generally, the notion of entanglement between two coproducts means that they verify
the entanglement equation}.
\end{defi}
\begin{theo}
Let $\D_1$ and $\D_2$ be two sub-spaces of a $k$-vector space $E$.
Suppose $(\D_1, \Delta_1, \delta_1)$ (resp. $(\D_1, \Delta_1, \hat{\delta}_1)$) and $(\D_2, \Delta_2, \delta_2)$ (resp. $(\D_2, \Delta_2, \hat{\delta}_2)$) are two codipterous coalgebras (resp. anti-codipterous coalgebras). The entanglement so obtained is still a codipterous coalgebra (resp. anti-codipterous coalgebra).
\end{theo}
\Proof
Denote by $\Delta_*$, the coassociative coproduct such that $\Delta_* := \Delta_1$ over $\D_1$, (resp. := $\Delta_2$ over $\D_2$) and $\delta_* := \delta_1$, (resp. $:= \delta_2$) over $\D_1$ (resp. over $\D_2$). The space 
$(\D := \D_1 + \D_2, \Delta_*, \delta_*)$ is a codipterous coalgebra. Similarly for the entanglement of two anti-codipterous coalgebras.
\eproof
\begin{defi}{[Pre-dendriform coalgebra]}
A $k$-vector space $\mathbb{D}$ equipped with three coproducts 
$\Delta, \delta, \hat{\delta}: \D \xrightarrow{} \D^{\otimes 2}$ verifying:
\begin{enumerate}
\item{Coas: $(\Delta \otimes id)\Delta = (id \otimes \Delta)\Delta$.}
\item{Codip: $(\Delta \otimes id)\delta = (id \otimes \delta)\delta$.}
\item{Anti-codip: $(id \otimes \Delta)\hat{\delta} = (\hat{\delta} \otimes id )\hat{\delta}$.}
\item{Entanglement equation: $(id \otimes \hat{\delta})\delta = (\delta \otimes id )\hat{\delta}$.}
\end{enumerate} 
is called a {\it{pre-dendriform coalgebra}}. Similarly, the coproducts $\delta$ and $\hat{\delta}$ are called {\it{bridges}}. Notice also that
if $\Delta = \delta + \hat{\delta}$, then the pre-dendriform coalgebra
becomes a dendriform coalgebra.
\end{defi}
\NB
With the notion of entanglement,
the reader may view the pre-dendriform coalgebra as the entanglement of the codipterous coalgebra $(\D, \Delta, \delta)$ with the anti-codipterous coalgebra $(\D, \Delta, \hat{\delta})$.
\NB
The axioms and terminology of dipterous algebras were discovered by J-L. Loday and 
M. Ronco \cite{LL,LR1}.
The axioms of pre-dendriform coalgebras were independently discovered by the author to describe the entanglement of two directed graphs. The denomination, pre-dendriform coalgebra, was also suggested by J-L. Loday \cite{LL}. If the coassociative law $\Delta$ is equal to $\delta + \hat{\delta}$, then
a pre-dendriform coalgebra turns out to be a dendriform coalgebra. In terms of directed graphs, this means that the directed graphs determined by the bridges $\delta$ and $\hat{\delta}$
realize a covering of the directed graph determined by $\Delta$. The covering becomes a tiling if the intersection of the arrow sets covering the graph determined by $\Delta$ is empty. The reader will notice then
the interest to break coassociativity in several coproducts. Fix an integer $n$, Let $(\delta_i)_{1 \leq i \leq n}$ be $n$ coproducts breaking the coassociativity of a coassociative coproduct $\Delta$, i.e. $n$ coproducts such that
$\sum_{1 \leq i \leq n} \ \delta_i = \Delta$. This decomposition can be viewed as a covering of the directed graph associated with $\Delta$ by
$n$ directed graphs associated with the coproducts $\delta_i$. The difficulty is to find convenient covering, i.e. convenient relations between coproducts such as by dualizing, the operads so obtained are Koszul. It is the case with the non-$\Sigma$ operad of dendriform algebras
and with the non-$\Sigma$ operad of dendriform trialgebras \footnote{See also J-L. Loday \cite{Lodaybh}.}. 

\begin{defi}{[Face, side]}
In this article, we will manipulate $k$-vector spaces \footnote{Such $k$-vector space will be called coassociative manifolds at the end of this paper.} $(M, (\delta_i)_{i \in I})$ equipped with coassociative coproducts $\delta_i$, $i \in I$. With the graphical viewpoint in mind, a {\it{face}} of such a space is $(M, (\delta_i)_{i \in J})$, $J \subset I$, a {\it{side}} is just a coassociative coalgebra $(M, \delta_{i_0})$, i.e. $J := \{ i_0 \}$.
\end{defi}

\noindent
A way to produce bridges
is to consider a channel map.
\begin{defi}{[Channel map]}
Let $E$ be a $k$-vector space and $C_1, C_2$ two subspaces of $E$, such that
$(C_1, \Delta_1, \epsilon_1)$ and $(C_2, \Delta_2, \epsilon_2)$, are
two coassociative coalgebras, resp. two bialgebras, resp. two Hopf algebras.
A linear map $\Phi$ is said to be a {\it{channel map}} if
$\Phi: \ C_1 \xrightarrow{} C_2 $ is invertible and
a coalgebra morphism, resp. a bialgebra morphism, resp. a Hopf algebra morphism, i.e. $\Delta_2 \Phi = (\Phi \otimes \Phi) \Delta_1$ and $\epsilon_2 \Phi = \epsilon_1$ and so on. If $C_1$ and $C_2$ are also algebras with units, then the channel must be unital.
\end{defi}
\subsection{Self-entanglement}
We now study the self-entanglement, i.e. the entanglement of two copies of a same codipterous coalgebra. A copy is produced thanks to a channel map. This kind of entanglement will yield associative dialgebras and associative trialgebras.
\begin{theo}
\label{dipt1}
Let $(C_1, \Delta_1, \epsilon_1)$ and $(C_2, \Delta_2, \epsilon_2)$
be two coassociative coalgebras of a $k$-vector space $E$, with $\Phi: \ C_1 \xrightarrow{} C_2 $, a channel map. Consider the subspace
$\D := C_1 + C_2$. Denote by $\Delta^*: \ \D \xrightarrow{} \D^{\otimes 2}$, the coproduct such that over $C_1$, $\Delta_* := \Delta_1$ 
and over $C_2$, $\Delta_* := \Delta_2$. Denote by $\delta_1: \ \D \xrightarrow{} \D^{\otimes 2}$, the coproduct such that over $C_1$, $\delta_1 := \Delta_1$ and $\delta_1 \Phi := (id \otimes \Phi)\Delta_1$. 
Then $(\D, \ \Delta_*, \ \delta_1)$ is a codipterous coalgebra. Moreover
$(\epsilon_1 \otimes id)\delta_1 =id$.
Denote by $\hat{\delta}_1: \ \D \xrightarrow{} \D^{\otimes 2}$ the coproduct such that over $C_1$, $\hat{\delta}_1 := \Delta_1$ and  $\hat{\delta}_1 \Phi := (\Phi \otimes id)\Delta_1$. 
Then $(\D, \ \Delta_*, \ \delta_1, \ \hat{\delta}_1)$ is a (chiral) pre-dendriform coalgebra.
Moreover $(id \otimes \epsilon_1)\hat{\delta}_1 =id$.
\end{theo}
\Proof
Notice that
$\Delta_*$ is coassociative. The equalities $(\epsilon_1 \otimes id)\delta_1 =id$ and $(id \otimes \epsilon_1)\hat{\delta}_1 =id$
are straightforward.
Let us check $(\Delta_* \otimes id) \delta_1 = ( id \otimes \delta_1) \delta_1$, which holds over $C_1$. Over $C_2$,
we observe that the right hand side is equal to 
$(id \otimes id \otimes \Phi)(\Delta_1 \otimes id)\Delta_1$
and the left hand side is equal to $(id \otimes id \otimes \Phi)(id \otimes \Delta_1)\Delta_1$.
We also easily obtain $(id \otimes \Delta_*) \hat{\delta}_1 = (\hat{\delta}_1 \otimes id ) \hat{\delta}_1$.
We can check that $(id \otimes \hat{\delta}_1) \delta_1 =
(\delta_1 \otimes id) \hat{\delta}_1$, holds over $C_1$. 
Over $C_2$, the equalities $(id \otimes \hat{\delta}_1) \delta_1 =
(id \otimes \Phi \otimes id)(id \otimes \Delta_1)\Delta_1$ and 
$(\delta_1 \otimes id) \hat{\delta}_1 = (id \otimes \Phi \otimes  id)(\Delta_1 \otimes id)\Delta_1$
prove that $\D$
is a (chiral) pre-dendriform coalgebra. 
\eproof
\Rk
\label{deuxponts}
To construct the bridges $\hat{\delta}_1$ and $\delta_1$ we decided to
produce a copy $C_2$ of $C_1$ via the channel map $\Phi$. Similarly, we
construct bridges from $C_2$ by using $\Phi^{-1}$
instead of $\Phi$. 
By reversing the point of view, we get
$\hat{\delta}_2 := \delta_2 : = \Delta_2$ over $C_2$. Over $C_1$,
$\hat{\delta}_2 \Phi^{-1}:= (\Phi^{-1} \otimes id)\Delta_2$ and
$\delta_2 \Phi^{-1}: = (id \otimes \Phi^{-1})\Delta_2$. Notice 
that results holding for the pre-dendriform coalgebra $(\D, \ \Delta_*, \ \delta_1, \ \hat{\delta}_1)$ are still valid with the pre-dendriform coalgebra $(\D, \ \Delta_*, \ \delta_2, \ \hat{\delta}_2)$.
\begin{theo}
With the hypotheses, the notation of the theorem \ref{dipt1}, and the
previous remark, the codipterous coalgebra $(\D, \ \Delta_*, \ \delta_1)$ is (chiral) entangled with the codipterous coalgebra $(\D, \ \Delta_*, \ \delta_2)$. Such an entanglement is called a self-entanglement, (since over $C_1$, $\delta_2 = \hat{\delta}_1 \Phi$) and is denoted by $[C_1 \underset{\delta_2}{{\stackrel{\delta_1}{\rightleftharpoons}}} C_1]$.
\end{theo}
\Proof
The theorem \ref{dipt1} and the previous remark imply that $(\D, \ \Delta_*, \ \delta_2, \ (\hat{\delta}_2))$ and $(\D, \ \Delta_*, \ \delta_1, \ (\hat{\delta}_1))$ are codipterous coalgebras.
The straightforward equality $(id \otimes \delta_2)\delta_1 = (\delta_1 \otimes id)\delta_2$ entails that $[C_1 \underset{\delta_2}
{{\stackrel{\delta_1}{\rightleftharpoons}}} C_1]$ is a self entangled codipterous coalgebra.
\eproof
\Rk
The two coalgebras $C_1$ and $C_2$ are called the boundaries of
the entangled codipterous coalgebras $[C_1 \underset{\delta_2}
{{\stackrel{\delta_1}{\rightleftharpoons}}} C_1]$. The boundary of a codipterous coalgebra will be denoted by $\partial$. Therefore, $\partial [C_1 \underset{\delta_2}
{{\stackrel{\delta_1}{\rightleftharpoons}}} C_1] = C_1 \cup C_2$.
\begin{prop}
With the hypotheses and the notation of the theorem \ref{dipt1}, we get
$(\hat{\delta}_1 +\delta_1) \Phi = (\Phi \otimes id)\Delta_1 + (id \otimes \Phi)\Delta_1$.
\end{prop}
\Proof
Straightforward.
\eproof
\NB \textbf{[Interpretation]}
Let us interpret the previous equality. Let $(C, \Delta)$ be
a coassociative coalgebra as well as an unital algebra with unit 1. We recall that a Leibniz coderivative is a linear map $D: C \xrightarrow{} C$
such that $D(1)=0$ and $\Delta D := (D \otimes id)\Delta + (id \otimes D)\Delta$. Suppose now that $(C, \delta, \hat{\delta} )$ is
a dendriform coalgebra. As $\Delta := \delta + \hat{\delta}$ is
coassociative, we get
$(\delta + \hat{\delta}) D := (D \otimes id)\Delta + (id \otimes D)\Delta$.
\begin{defi}{[Leibniz coderivative on a pre-dendriform coalgebra]}
With the hypotheses and the notation of the theorem \ref{dipt1}, suppose $(\D, \Delta_*,\delta_1, \hat{\delta}_1 )$ is a pre-dendriform coalgebra as well as an algebra with unit $1$. A Leibniz coderivative on the pre-dendriform coalgebra $\D$
is a linear map $D_I: C_1 \xrightarrow{} C_2$ such that 
$(\delta_1 + \hat{\delta}_1) D_I := (D_I \otimes id)\Delta_1 + (id \otimes D_I)\Delta_1$ and $D_I(1)=0$.
\end{defi}
\begin{prop}
With the hypotheses and the notation of the theorem \ref{dipt1}, the map $D_I := \Phi - id$ is a Leibniz coderivative on the pre-dendriform coalgebra $\D$. If $\D$ is a pre-dendriform bialgebra then $D_I$ is an Ito derivative.
\end{prop}
\Proof
Let $(\D, \Delta_*,\delta_1, \hat{\delta}_1 )$ be the pre-dendriform coalgebra described in theorem \ref{dipt1}. By definition of the channel map, $D_I(1)=0$ is straightforward. The Leibniz coderivative equality comes from the definitions of the bridges. If $(\D, \Delta_*,\delta_1, \hat{\delta}_1 )$ is a pre-dendriform bialgebra, as $\Phi$ is a homomorphism,
we easily check the Ito property: $D_I(xy) - D_I(x)D_I(y) = xD_I(y) + D_I(x) y$, for all
$x,y \in C_1$.
\eproof
\begin{exam}{[Axioms of codipterous coalgebras: a graphical interpretation ]}
Let $E$ be the $k$-vector space generated by $a,b,c,d,x,y,z$ and $u$.
Let $\F$ be the subspace generated by $a,b,c$ and $d$. Define
$\Delta: \F \xrightarrow{} \F^{\otimes 2}$ such that:
$
\Delta a = a \otimes a + b \otimes c, \ \
\Delta b = a \otimes b + b \otimes d, \ \
\Delta c = d \otimes c + c \otimes a, \ \
\Delta d = d \otimes d + c \otimes b.
$
\begin{center}
\includegraphics*[width=3.5cm]{graph2.eps}

\textbf{\begin{scriptsize} Geometric support of $\F$. \end{scriptsize}}
\end{center}
Consider the coassociative coalgebra $\F'$ generated by $x,y,z$ and $u$, with coproduct $\Delta': \F' \xrightarrow{} \F'^{\otimes 2}$ defined via the channel map
$\Phi: \F \xrightarrow{} \F'$, i.e. $\Delta' \Phi : = (\Phi \otimes \Phi)\Delta$, with $\Phi(a) = x$, $\Phi(b) = y$, $\Phi(c) = z$, $\Phi(d) = u$.
The geometric support of the entangled codipterous coalgebras constructed in the theorem \ref{dipt1} has two directed graphs, $\F$ and $\F'$, with four bridges
between them $\delta_1$, $\hat{\delta}_1$, $\delta_2$ and $\hat{\delta}_2$.
The bridge $\delta_1$ is defined over $\F'$ such that:
$
\delta_1 x = a \otimes x + b \otimes z, \ \
\delta_1 y = a \otimes y + b \otimes u, \ \
\delta_1 z = d \otimes z + c \otimes x, \ \
\delta_1 u = d \otimes u + c \otimes y.
$
The bridge $\hat{\delta}_1$ is defined over $\F'$ such that:
$
\hat{\delta}_1 x = x \otimes a + y \otimes c, \ \
\hat{\delta}_1 y = x \otimes b + y \otimes d, \ \
\hat{\delta}_1 z = u \otimes c + z \otimes a, \ \
\hat{\delta}_1 u = u \otimes d + z \otimes b.
$
We get $\partial [\F \underset{\delta_2}{{\stackrel{\delta_1}{\rightleftharpoons}}}\F] 
= \F \cup \F'$.
The geometric support of  $[\F \underset{\delta_2}{{\stackrel{\delta_1}{\rightleftharpoons}}}\F] 
$ is:
\begin{center}
\includegraphics*[width=7cm]{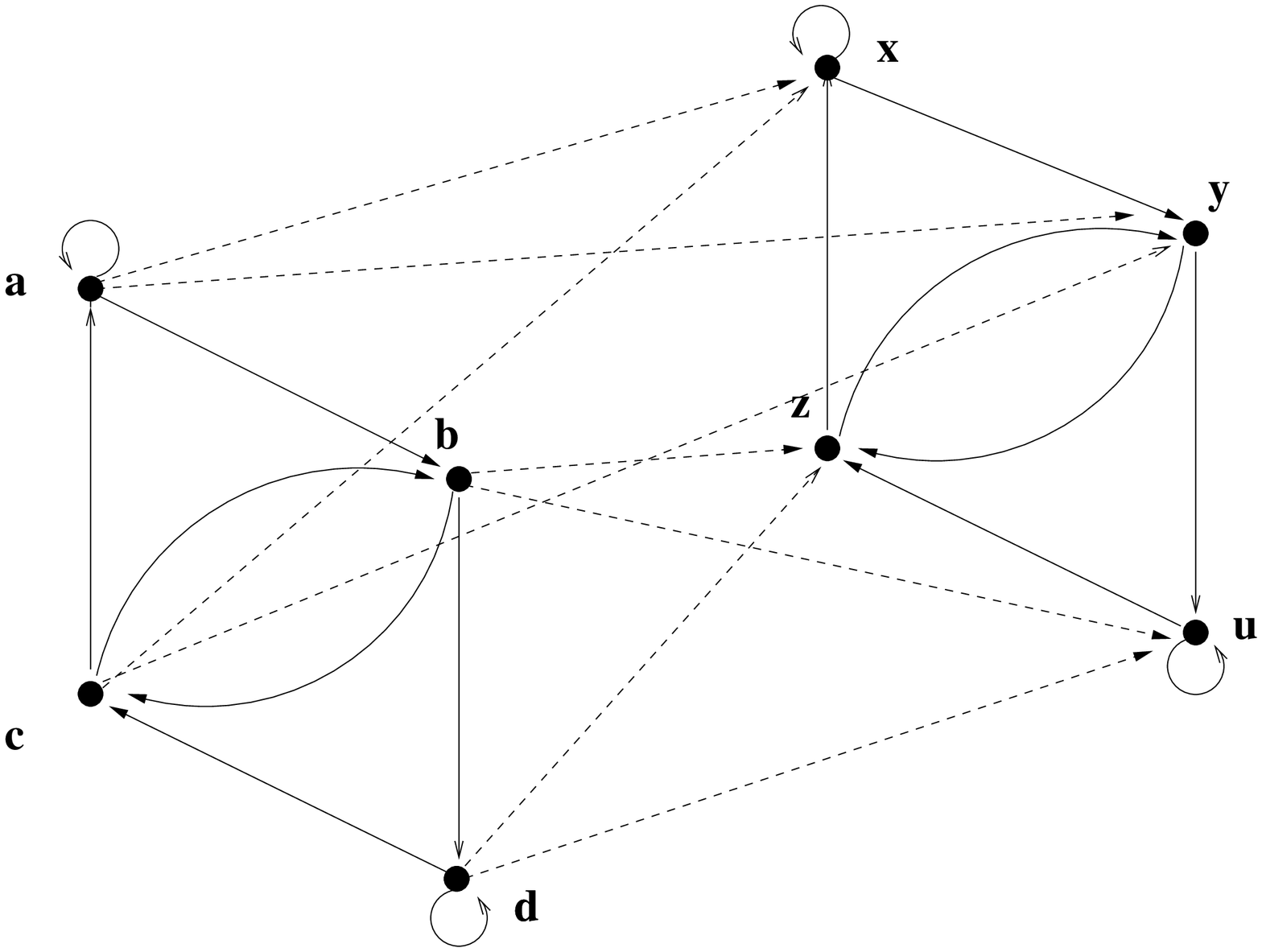}

\textbf{\begin{scriptsize} Geometric support of $[\F \underset{\delta_2}{{\stackrel{\delta_1}{\rightleftharpoons}}}\F] 
$. Among the bridges, only
$\delta_1$ is represented.\end{scriptsize}}
\end{center}
Since the chanel map has no fixed point, as a $k$-vector space, $[\F \underset{\delta_2}{{\stackrel{\delta_1}{\rightleftharpoons}}}\F] 
:= F \oplus \Phi(F) :=E$.
\end{exam}
If the channel has fixed points, we have to be careful.
\begin{prop}
Let $E':= C_1 + C_2$ be a $k$-vector space, where
$C_1$ is a coassociative coalgebra generated by $(x_{i})_{i \in I}$ as a subspace, with coproduct $\Delta x_{i} = \sum x_{i, (1)}
\otimes x_{i, (2)}$ and let $\Phi':C_1 \xrightarrow{} C_2$ be a channel. Define $\Delta_*:=\Delta$ over $C_1$ and $\Delta_*\Phi':=(\Phi' \otimes \Phi')\Delta$ over $C_2$ and $\delta := \Delta$ over $C_1$ verifying $\delta \Phi' := (id \otimes \Phi')\Delta$. Fix $i\in I$. Suppose $x_{i}$ is fixed by the channel
map $\Phi'$. 
Then the sub-coalgebra generated by $x_{i}$
is fixed by $\Phi'$.
\end{prop}
\Proof
Fix $i \in I$ and suppose $x_{i}$ is fixed by the channel
map $\Phi'$. Then $\Delta(x_{i}):= \sum x_{i,  (1)} \otimes x_{l,  (2)}$ has to be equal to $\delta x_{i}:= \sum x_{i, (1)} \otimes \Phi'(x_{i, (2)})$ which entails that $\Phi'((x_{i, (2)}))$ is fixed by the channel. Similarly, $\hat{\delta} x_{i}:= \Delta(x_{i})$ implies that $\Phi'((x_{i, (1)}))$ is fixed by $\Phi'$. Denote by $S$, the coalgebra generated by $x_i$ and by $S_1$ the subspace of fixed points of $\Phi'$ belonging to $S$. We notice that $x_i \in S_1$, $S_1 \subset S$ and $\Delta S_1 \subset S_1 \otimes S_1$. Therefore $S_1$ is a sub-coalgebra containing $x_i$ thus $S_1 =S$.
\eproof
\begin{theo}
With the hypotheses and the notation of the theorem \ref{dipt1}, the pre-dendriform coalgebra
$(\D, \delta_1, \hat{\delta}_1)$ is a coassociative codialgebra
and $(\D, \Delta_*, \delta_1, \hat{\delta}_1)$ is a coassociative cotrialgebra.
\end{theo}
\Proof
Let us prove that the pre-dendriform coalgebra
$(\D_1, \delta_1, \hat{\delta}_1)$ is a coassociative co-dialgebra. The bridges
$\delta_1$ and $\hat{\delta}_1$ are coassociative. This holds over
$C_1$. Over $C_2$, since $\delta_1 = \Delta_1$ over $C_1$, we get
$(\delta_1 \otimes id)\delta_1 = (\Delta_1 \otimes id)\delta_1 =
(id \otimes \delta_1)\delta_1$, which proves the coassociativity of
$\delta_1$. Similarly, by using the anti-codipterous coalgebra axiom, we get the coassociativity of $\hat{\delta}_1$. The entanglement equation
$(id \otimes \hat{\delta}_1) \delta_1 =
(\delta_1 \otimes id) \hat{\delta}_1$ is assumed in the definition of a pre-dendriform coalgebra.  Over $C_2$, the equalities
$(\delta_1 \otimes id) \delta_1 = (\hat{\delta}_1 \otimes id)\delta_1$ 
and $(id \otimes \delta_1) \hat{\delta}_1 = (id \otimes \hat{\delta}_1)\hat{\delta}_1$ 
hold, since $\delta_1 = \hat{\delta}_1 = \Delta_1$ over $C_1$.
This proves that $(\D, \delta_1, \hat{\delta}_1)$ is a coassociative co-dialgebra. Similarly, the proof that $(\D, \Delta_*, \delta_1, \hat{\delta}_1)$ is a coassociative co-trialgebra is complete by checking the axioms of the definition \ref{cotri}.
\eproof
\Rk
The previous theorem is still valid for the pre-dendriform coalgebra $(\D, \Delta_*, \delta_2, \hat{\delta}_2)$.
\begin{coro}
With the hypotheses and the notation of theorem \ref{dipt1},
let $(A,m)$ be an associative algebra over $k$. Denote by $L(\D, A)$,
the $k$-vector space of the linear maps which map $\D$ into $A$.
Let $f,g \in L(\D, A)$. Denote by $f \perp g :=m(f \otimes g)\Delta^*$, $f \dashv g := m(f \otimes g) \hat{\delta}_1$ and 
$f \vdash g := m(f \otimes g) \delta_1$, the convolution products.
Then $(L(\D, A),\dashv,\vdash)$
is an associative dialgebra. Define for all
$f,g \in L(\D, A)$, the bracket $[f,g] := f \dashv g  - f \vdash g $.
This bracket embeds $L(\D, A)$ into a Leibniz algebra.
Moreover $(L(\D, A),\dashv,\vdash,\perp)$ is an associative trialgebra.
\end{coro}
\Proof
Straightforward.
\eproof
\NB
With the hypotheses and the notation of theorem \ref{dipt1}, notice
there exist directed graphs, which are geometric supports of
particular pre-dendriform coalgebras $(\D, \Delta_*, \delta_1, \hat{\delta}_1)$, constructed from coassociative coalgebras
via a channel map $\Phi$, whose algebra of linear maps
$L(\D, A)$, where $A$ is an algebra with unit $1_A$, has a structure of associative di or trialgebra.
The Leibniz bracket allows us to deal with differential structures on 
these directed graphs. Notice that on the boundary $C_1$ of such a pre-dendriform
coalgebra, the coproducts $\delta_1 = \hat{\delta}_1 = \Delta_1$, therefore
the Leibniz algebra turns out to be a Lie algebra on the boundary $C_1$. With associative trialgebras, we can also construct {\it{``non commutative version'' of Poisson algebras}}. Defined in \cite{LodayRonco},
a non commutative version of Poisson algebra $P$ is a $k$-vector space equipped 
with a Leibniz bracket and an associative operation $x \bullet y$ (not necessarily commutative) such that its relationship with the Leibniz bracket is given by:
$$ \forall x,y,z \in P, \ [x \bullet y,z] = x \bullet [y,z] + [x,z] \bullet y, \ \ [x, y \bullet z -z \bullet y] = [x, [y,z]].$$
By defining for all $f,g \in L(\D, A)$, $f \bullet g := f \perp g$, the associative trialgebra $L(\D, A)$ becomes a Poisson algebra. 
\Rk
In \cite{Loday}, J-L. Loday defines the notion of a bar-unit in an associative dialgebra $(X, \vdash, \dashv)$. An element $e \in X$ is said to be a {\it{bar-unit}} of $X$, if for all $x \in X$, $x \dashv e =x = e \vdash x$. The set of bar-units is called the {\it{halo}}. Let us
still work with the hypotheses and the notation of the theorem \ref{dipt1}. Denote by $\eta_A: k \xrightarrow{} A$ such that $\lambda \mapsto \lambda 1_A$.
If the involved coassociative coalgebra $C_1$ has a counit $\epsilon: \ C_1 \xrightarrow{} k$, then the map $e:= \eta_A \epsilon_*$, where $\epsilon_* = \epsilon$ over $C_1$ and $\epsilon_*\Phi = \epsilon$ over $C_2$, is a bar-unit in $L(D,A)$, since if $f \in L(D,A)$, then $e \vdash f  := m(e \otimes f)\delta_1 = m(id \otimes f)( \eta_A \epsilon_* \otimes id) \delta_1 = f$. Similarly, $f \dashv v  := m(f \otimes e)\hat{\delta}_1 = m(f \otimes id)(  id \otimes \eta_A \epsilon_*) \hat{\delta}_1 = f$. As the counit $\epsilon$ is unique, there will exist an unique bar-unit on $L((\D, \Delta_*, \delta_1, \hat{\delta}_1), A)$, i.e. the halo is a singleton \footnote{By entangling several codialgebras, constructed with this method, we could obtain a larger halo.}.
\begin{theo}
Consider the hypotheses and the notation of the theorem \ref{dipt1}.
If $C_1$, (resp. $C_2$) is generated as a $k$-vector space by a basis $(v_i)_{1 \leq i \leq \dim C_1}$, (resp. $(w_i)_{1 \leq i \leq \dim C_1}$), with $\dim C_1 < \infty$, we denote 
by $(v^*_i)_{1 \leq i \leq \dim C_1}$, (resp. $(w^*_i)_{1 \leq i \leq \dim C_1}$), the dual basis, i.e. $v^*_i(v_j)= 1$ (resp. $w^*_i(w_j)= 1$)
if $i=j$ and 0 otherwise. 
Denote by $B_{ij}^k$ and $C_{ij}^k$ the structure constants, i.e.
scalars such that $[ w^*_i,v^*_j ] = \sum_k B_{ij}^k w^*_k $
and $[ v^*_i, v^*_j ]= \sum_k C_{ij}^k v^*_k $.
The structure constants are determined by the geometric support of
the pre-dendriform coalgebra $\D$, i.e. $ v^*_i \otimes v^*_j \not=0$ if and only if
$v_i \xrightarrow{} v_j$ is an arrow of the directed graph associated with $\D$, similarly for $ w^*_i \otimes v^*_j$ and $ v^*_i \otimes w^*_j$, for all $1 \leq i,j \leq \dim C_1$. Moreover, $[v^*_j,w^*_i] =0$, for all $1 \leq i,j \leq \dim C_1$.
\end{theo}
\Proof
Fix $i,j \in 1 \leq i \leq \dim C_1$. Denote by $m$, the product of the
field $k$.
Let $v_j$, (resp. $w_i$) be an element of the basis of $C_1$ (resp. of $C_2$). By definition, $[v^*_i,w^*_j] := v^*_i \dashv w^*_j - w^*_j \vdash v^*_i = m (v^*_i \otimes w^*_j)\hat{\delta}_1 - m (w^*_j \otimes v^*_i)\delta_1$, hence  $[v^*_i,w^*_j] = 0$.
\eproof
\begin{exam}{}
In the case of the pre-dendriform coalgebra $\D := [\F \dpe \F]$,
on the boundary $C_1$, the Leibniz bracket turns out to be the Lie bracket determined by the convolution product associated with $\Delta_* := \Delta_1$. For instance we get 
$[a^*,b^*] = b^*$, $[b^*,c^*] = a^* - d^*$, which behaves as a Lie bracket and  $[ a^*,x^*] = 0$, $[ x^*,a^*] = x^*$, $[ y^*,c^* ] = x^*-u^*$,
$[c^*,y^*]= 0$, out of the boundary $C_1$, which behaves as a Leibniz bracket.
\end{exam}
\begin{theo}
Consider the hypotheses and the notation of the theorem \ref{dipt1}. In addition, suppose that the coalgebras involved are bialgebras.
As the channel $\Phi$ is a homomorphism, then so are the bridges
$\delta_1$ and $\hat{\delta}_1$.
\end{theo}
\Proof
Let us prove the theorem for the bridge $\delta_1$.
Over $C_1$, $\delta_1 := \Delta_1$, is obviously an unital homomorphism.
Let $x,y \in C_2$. There exist unique $a,b \in C_1$, such that
$\Phi(a) =x$ and $\Phi(b) =y$.
Then $\delta_1(xy) = \delta_1\Phi(ab) := (id \otimes \Phi)\Delta_1(ab)$.
With the Sweedler's notation, we write $\Delta_1(a):= \sum a_{(1)} \otimes a_{(2)}$
and $\Delta_1(b):= \sum b_{(1)} \otimes b_{(2)}$.
Therefore, $(id \otimes \Phi)\Delta_1(ab) := \sum a_{(1)}b_{(1)} \otimes \Phi(a_{(2)}b_{(2)}) = (\delta_1 \Phi(a))(\delta_1 \Phi(b))$, 
which proves the homomorphism property of the bridge $\delta_1$. Since the channel $\Phi$ is by definition unital, so is the bridge $\delta_1$.
\eproof
\begin{theo}
\label{hop}
Consider the hypotheses and the notation of the theorem \ref{dipt1}. In addition, suppose the coalgebras involved are Hopf algebras
with antipodes $S_1$ and $S_2$
and recall that the channel $\Phi$ is a morphism of Hopf algebra.
Denote by $S_*$, the linear map that equals to $S_1$ over $C_1$
and to $S_2$ over $C_2$. 
Denote by $\sigma_1$, the linear map that equals to $S_1$ over $C_1$ and to
$\Phi^{-1} S_2$ over $C_2$. Then
the bridges $\delta_1$ and $\hat{\delta}_1$ verifies
$m(id \otimes \sigma_1)\delta_1 = 1\epsilon_*$
and $m(\sigma_1 \otimes id)\hat{\delta}_1 = 1\epsilon_*$. Similarly,
denote by $\sigma_2$, the linear map that equals to $S_2$ over $C_2$ and to
$\Phi S_1$ over $C_1$. Then
the bridges $\delta_2$ and $\hat{\delta}_2$ verify
$m(id \otimes \sigma_2)\delta_2 = 1\epsilon_*$
and $m(\sigma_2 \otimes id)\hat{\delta}_2 = 1\epsilon_*$. Therefore,
the pre-dendriform bialgebras $(\D, \Delta_*,\delta_1, \hat{\delta}_1 )$, $(\D, \Delta_*,\delta_2, \hat{\delta}_2 )$ and $(\D, \Delta_*,\delta_1, \hat{\delta}_2 )$ are $L$-Hopf algebras.
\end{theo}
\Proof
Each side of the pre-dendriform bialgebra $(\D, \delta_1, \hat{\delta}_1 )$, i.e. $(\D, \delta_1)$ and $(\D, \hat{\delta}_1)$, has a map $\sigma_1$, 
defined in this theorem, such that $id \vdash \sigma_1= e$ and $\sigma_1 \dashv id =e$, where $e$ is the bar unit of $L(\D, \D)$. As their coproducts verify the entanglement equation, $(\D, \delta_1, \hat{\delta}_1 )$ can be viewed as a $L$-Hopf algebra \cite{Coa}. Similarly for the space $(\D, \delta_2, \hat{\delta}_2 )$. 
For the space $(\D, \Delta_*,\delta_1, \hat{\delta}_2 )$, as the entanglement equation $ (id \otimes \hat{\delta}_2) \delta_1 = (\delta_1 \otimes id) \hat{\delta}_2$ is verifyied, the space
$(\D, \Delta_*,\delta_1, \hat{\delta}_2 )$ has the structure of $L$-Hopf algebra \footnote{It is not a coassociative codialgebra.}.
\eproof

\noindent
Let $(\D, \Delta_*, \delta_1, \hat{\delta}_1)$ be the pre-dendriform coalgebra defined in the theorem \ref{dipt1}.
Inspired with results in \cite{Hudson1}, we define 
$\overrightarrow{d}: \ \D \xrightarrow{} \D \otimes \D$, such that $\overrightarrow{d}
=\Delta_*  - \delta_1$ and
$\overleftarrow{d}: \ \D \xrightarrow{} \D \otimes \D$, such that $\overleftarrow{d}
=    \Delta_*  - \hat{\delta}_1. $
\begin{theo}
Let $(\D, \delta_1, \hat{\delta}_1)$ be the pre-dendriform coalgebra
defined in theorem \ref{dipt1}. 
We get:
$(id \otimes \overleftarrow{d})\overrightarrow{d}= (\overrightarrow{d} \otimes id)\overleftarrow{d}.$ Moreover, $(\overrightarrow{d} \otimes id) \Delta_* = (id \otimes \Delta_*)\overrightarrow{d}$ and $(id \otimes \overleftarrow{d} ) \Delta_* = (\Delta_* \otimes id)\overleftarrow{d}$.
\end{theo}
\Proof
Straightforward by using the following equality
$(id \otimes \hat{\delta})\delta = (\delta \otimes id )\hat{\delta}.$
\eproof
\begin{defi}{[bimodule]}
Let $(\D, \Delta_*, \delta_1, \hat{\delta}_1)$ be the pre-dendriform coalgebra
defined in the theorem \ref{dipt1}. Recall that the coproducts $\delta_1$ and $\hat{\delta}_1$ are coassociative. Suppose $\D$ is a pre-dendriform bialgebra.
Embed $\D ^{\otimes 2}$
into a {\it{$\D$-bimodule}} by defining for all
$c, x, y \in \D$,
\begin{eqnarray*}
x \ \hat{\circ} \ \overleftarrow{d}(c) &=& \hat{\delta}_1(x)\overleftarrow{d}(c); \ \ \ \ \overleftarrow{d}(c) \ \hat{\circ} \ y = \overleftarrow{d}(c) \hat{\delta}_1(y).
\end{eqnarray*}
Similarly, $\D ^{\otimes 2}$ can be embedded into another {\it{$\D$-bimodule}} by defining for all
$c, x, y \in \D$,
\begin{eqnarray*}
x \circ \overrightarrow{d}(c) &=& \delta_1(x)\overrightarrow{d}(c); \ \ \ \ \overrightarrow{d}(c)\circ y = \overrightarrow{d}(c)\delta_1(y).
\end{eqnarray*}
\end{defi}
\begin{theo}
\label{dipt3}
Let $(\D, \delta_1, \hat{\delta}_1)$ be the pre-dendriform bialgebra
defined in theorem \ref{dipt1}. 
Then $\overleftarrow{d},\overrightarrow{d}$ are Ito derivatives.
\end{theo}
\Proof
Let $x,y \in \D$.
We have $\overleftarrow{d}(1) =0= \overrightarrow{d}(1)$. Moreover,
$\overrightarrow{d}(x)\overrightarrow{d}(y) = \Delta_*(xy) + \delta_1(xy) - \Delta_*(x)\delta_1(y) - \delta_1(x)\Delta_*(y)$, i.e.
$ \overrightarrow{d}(xy)= \overrightarrow{d}(x)\overrightarrow{d}(y) + \overrightarrow{d}(x) \ \circ \ y + x \ \circ \ \overrightarrow{d}(y).$
Similarly,
$ \overleftarrow{d}(xy)= \overleftarrow{d}(x)\overleftarrow{d}(y) + \overleftarrow{d}(x)\hat{\delta}_1(y) + \hat{\delta}_1(x)\overleftarrow{d}(y)$, that is $ \overleftarrow{d}(xy)= \overleftarrow{d}(x)\overleftarrow{d}(y) + \overleftarrow{d}(x) \ \hat{\circ} \ y + x \ \hat{\circ} \  \overleftarrow{d}(y)$.
\eproof
\Rk
Let $(\D, \delta_1, \hat{\delta}_1)$ be the pre-dendriform coalgebra
defined in theorem \ref{dipt1}.
If $x \in \D$ verifies $\Delta_*(x) = x \otimes 1 + 1 \otimes x$,
then $( \overleftarrow{d} \otimes id)\overrightarrow{d}(x)=0$.
\subsection{Self-tilings and dendriform coalgebras}
We now obtain dendriform coalgebras from the self-entanglement and apply
these results on particular coalgebras.
\begin{theo}
With the hypotheses, the notation of the theorem \ref{dipt1},
as usual,
define the bridge $\delta_d :=\Delta$ over $C_1$ such that
$\delta_d\Phi := (id \otimes \Phi)\Delta$ over $C_2$. 
Define also $\hat{\delta}_d \Phi:=(\Phi \otimes id)\Delta$ over $C_2$ and $\hat{\delta}_d := 0$ otherwise. Define also $\Delta_*\Phi:
=(\Phi \otimes \Phi) \Delta$ over $C_2$ and $\Delta_*:
= \Delta $ over $C_1$.
Then the coproduct $\bar{\Delta} := \delta_d + \hat{\delta}_d$ is coassociative and equipped with the coproducts $\delta_d$,  $\hat{\delta}_d$, $E$ becomes a dendriform coalgebra.
\end{theo}
\Proof
Straightforward by checking axioms.
\eproof
\Rk
The dendriform coalgebra is trivial on the boundary $C_1$. The directed graph constructed from $\bar{\Delta}$ is tiled
by the directed graphs constructed from its decomposition, that is by the coproducts $\delta_d$ and $\hat{\delta}_d$.

\noindent
Let us yield an application of entanglement given by two channels, inspired from a work of C. Cibils \cite{Cibils}. Fix an integer $n>0$ and an invertible scalar $q \in k$.
Let $E$ be a $k$-vector space generated by $(a_i)_{0 \leq i \leq n-1}$ and by $(x_i)_{0 \leq j \leq n-1}$. Define the following coassociative coproduct $\Delta$ on $k \bra (a_i)_{0 \leq i \leq n-1} \ket$ such that
$\Delta(a_i) = \sum_{j+k=i} a_j \otimes a_k$. 
Define the channel map $\Phi: k \bra (a_i)_{0 \leq i \leq n-1} \ket \xrightarrow{} k \bra (x_i)_{0 \leq i \leq n-1} \ket$ such that $\Phi(a_i)= q^{-i}
x_i.$ Then the bridge $\hat{\delta}$ verifies
$\hat{\delta}(x_i) = \sum_{j+k=i} q^k x_j \otimes a_k$ and $\hat{\delta}(a_i) = \sum_{j+k=i} a_j \otimes a_k$ and
is coassociative. The map $\epsilon(a_i)= 0$, if $i \not= 0$, $\epsilon(a_0)= 1$ and $\epsilon(x_i)= 0$ is a right counit. Moreover $(E, \Delta_*, \hat{\delta})$, where $\Delta_*\Phi := (\Phi \otimes \Phi)\Delta$ over $k \bra (x_i)_{0 \leq i \leq n-1} \ket$ is an anti-codipterous coalgebra. Similarly,
define the channel map $\Psi: k \bra (a_i)_{0 \leq i \leq n-1} \ket \xrightarrow{} k \bra (x_i)_{0 \leq i \leq n-1} \ket$ such that $\Psi(a_i)= 
x_i.$ Then the bridge $\delta$ verifies
$\delta(x_i) = \sum_{j+k=i} a_j \otimes x_k$ and $\delta(a_i) = \sum_{j+k=i} a_j \otimes a_k$ and
is coassociative. The map $\epsilon(a_i)= 0$, if $i \not= 0$, $\epsilon(a_0)= 1$ and $\epsilon(x_i)= 0$ is a left counit. Moreover $(E, \Delta^*, \delta)$, where $\Delta^*\Psi := (\Psi \otimes \Psi)\Delta$ over $k \bra (x_i)_{0 \leq i \leq n-1} \ket$ is a codipterous coalgebra.
\begin{prop}
The $k$-vector space $E$ equipped with the two
coproducts $\delta$, $\hat{\delta}$
is a coassociative co-dialgebra and $(E, \Delta_*, \hat{\delta})$ and $(E, \Delta^*, \delta)$ are
chiral entangled.
\end{prop}
\Proof
Notice \footnote{This can be considered as a generalisation of the self-entanglement and works thanks to the particular form of the involved coproduct.} that $\Delta^*:=\Delta_*$ and that $(id \otimes \hat{\delta}) \delta = (\delta  \otimes id)\hat{\delta} $. The proof is complete by checking the axioms of a co-dialgebra. 
\eproof
\begin{theo}
Define now the coproduct $\delta$ such that for
all $0 \leq i \leq n-1$, $\delta(a_i) :=\Delta(a_i)$ and 
$\delta(x_i) := \sum_{j+k=i} a_j \otimes x_k$. Similarly, define the coproduct $\hat{\delta}$ such that for
all $0 \leq i \leq n-1$, $\hat{\delta}(a_i) :=0$ and 
$\hat{\delta}(x_i) := \sum_{j+k=i} q^k x_j \otimes a_k$. Then $(E, \delta, \hat{\delta})$ is a dendriform coalgebra.
\end{theo}
\Proof
Notice that $\Delta := \hat{\delta} + \delta$ is coassociative \cite{Cibils}. It is trivial on $(a_k)_{0 \leq k \leq n-1}$. Let us check the axioms on $(x_k)_{0 \leq k \leq n-1}$.
For the codipterous coalgebra axiom: $(id \otimes \delta)\delta(x_i) :=
\sum_{j+(l+m)=i} a_j \otimes a_l \otimes x_m$ which is equal to $(\Delta \otimes id)\delta(x_i) := \sum_{l'+j'+k=i} a_{j'} \otimes a_{l'} \otimes x_k$ since the sum is over all possible decompositions of $i$ into three integers. Similarly for the anti-codipterous coalgebra axiom: $(\hat{\delta} \otimes id)\hat{\delta}(x_i) :=
\sum_{j' + l'+k=i} q^k (q^{l'} x_{j'} \otimes a_{l'})\otimes a_k$
which is equal to $(id \otimes \Delta)\hat{\delta}(x_i) := \sum_{j+(h'+m')=i} q^{(h'+m')} x_j \otimes a_{h'} \otimes a_{m'}$, since the sum is over all possible decompositions of $i$ into three integers.
The entanglement equation is also verified. 
\eproof
\section{Entanglement of two different codipterous coalgebras}
Another way to construct entangled codipterous coalgebra is to start, for instance, from an achiral coassociative $L$-coalgebra, i.e. a $L$-coalgebra whose two coproducts $\Delta$ and $\tilde{\Delta}$ verify the entanglement equation whatever the position of the two coproducts are, i.e. $(\Delta \otimes id)\tilde{\Delta}= (id \otimes \tilde{\Delta})\Delta$ and $(id \otimes \Delta)\tilde{\Delta}= (\tilde{\Delta} \otimes id)\Delta$.
\begin{theo}
\label{dipt2}
Let $G$ be a sub-vector space of a $k$-vector space $E$, with $\dim E \geq
2 \dim G$, if $\dim G$ is finite.
Suppose $(G, \Delta, \tilde{\Delta}, \epsilon, \tilde{\epsilon})$ is an achiral coassociative $L$-coalgebra.
Define a channel map $\Phi: G \xrightarrow{} E$, verifying $\tilde{\Delta} \Phi := (\Phi \otimes \Phi) \tilde{\Delta}$. Define $C_1 :=G$, $C_2 := \Phi(G)$, $\Delta_1 :=\Delta$ and
$\tilde{\Delta}_2 := \tilde{\Delta} \Phi$. Define $\D := C_1 + C_2$ and 
$\delta_1: \D \xrightarrow{} \D^{\otimes 2}$ such that $\delta_1:= \Delta_1$ over
$C_1$ and $\delta_1 \Phi = (id \otimes \Phi) \Delta_1$ over $C_2$. Similarly, define 
$\tilde{\delta}_2: \D \xrightarrow{} \D^{\otimes 2}$, verifying $\tilde{\delta}_2= \tilde{\Delta}_2$ over
$C_2$ and $\tilde{\delta}_2 \Phi^{-1} = (\Phi^{-1} \otimes id) \tilde{\Delta}_2$ over $C_1$.
Denote by $\Delta_* := \Delta_1$ over $C_1$ and $\Delta_* = \tilde{\Delta}_2$ over $C_2$.
The space $(\D,\Delta_*, \delta_1, \tilde{\delta}_2)$ is the entanglement of two codipterous coalgebras $(\D,\Delta_*, \delta_1)$ and $(\D,\Delta_*, \tilde{\delta}_2)$. This (chiral) entanglement is denoted by $[G \dpi \tilde{G}]$, where $\tilde{G} = \Phi(G)$.
\end{theo}
\Proof
We check that $(\Delta_1 \otimes id)\delta_1 =(\Delta_* \otimes id)\delta_1 = (id \otimes \delta_1)\delta_1$ and 
$(\tilde{\Delta}_2 \otimes id)\tilde{\delta}_2 =(\tilde{\Delta}_* \otimes id)\tilde{\delta}_2 = (id \otimes \tilde{\delta}_2)\tilde{\delta}_2$. Moreover $(\tilde{\delta}_2 \otimes id)\delta_1 = (id \otimes \delta_1)\tilde{\delta}_2$ is straightforward.
\eproof
\begin{theo}
Consider the hypotheses and the notation of theorem \ref{dipt2}. In addition, suppose the coalgebras involved are Hopf algebras
with antipodes $S_1$ and $\tilde{S}_2$
and recall that the channel $\Phi$ is a morphism of Hopf algebras.
Denote by $S_*$, the linear map equals to $S_1$ over $C_1$
and to $\tilde{S}_2$ over $C_2$. 
Denote by $\sigma_1$, the linear map equals to $S_1$ over $C_1$ and to
$\Phi^{-1} \circ \tilde{S}_2$ over $C_2$. Define the map $\epsilon_*:= \epsilon$ over $C_1$ and $\epsilon_*\Phi:= \tilde{\epsilon}$ over $C_2$, 
the bridge $\delta_1$ verifies
$m(id \otimes \sigma_1)\delta_1 = 1\epsilon_*$. Similarly,
denote by $\tilde{\sigma}_2$, the linear map equals to $\tilde{S}_2$ over $C_2$ and to
$\Phi \circ S_1$ over $C_1$. Then
the bridge $\tilde{\delta}_2$ verifies
$m(id \otimes \tilde{\sigma}_2)\tilde{\delta}_2 = 1\epsilon_*$. Therefore,
the space $[G \dpi \tilde{G}]$ is a $L$-Hopf algebra.
\end{theo}
\Proof
The proof is similar to the proof of the theorem \ref{hop}. Even if such an entanglement does not yield a coassociative trialgebra, the element $e' :=\eta_A \epsilon_*$, with $\epsilon_*$ is defined in this theorem, plays the r\^ole of a bar-unit.
\eproof
\begin{theo}
\label{dipt4}
Let $(\D,\Delta_*, \delta_1, \tilde{\delta}_2)$ be the entangled codipterous coalgebra defined in theorem \ref{dipt2}. Let us define the coproduct $\Delta^*$ (resp. $\tilde{\Delta}^*$) such that 
$\Delta^* := \Delta_1$ (resp. $:= \tilde{\Delta}_1$) over $C_1$ and equals to $\Delta_2$ such that $\Delta_2 \Phi:= (\Phi \otimes \Phi) \Delta_1$ (resp. equals to $\tilde{\Delta}_2$ such that $:= \tilde{\Delta}_2\Phi:= (\Phi \otimes \Phi) \tilde{\Delta}_1$) over $C_2$. The coalgebra
$(C_1 + C_2,\Delta^*, \tilde{\Delta}^*)$ is then an achiral $L$-coalgebra. Recall that the bridges $\delta_1$ and $\hat{\delta}_2$
are defined by $\delta_1 := \Delta_1$ over $C_1$, such that $\delta_1 \Phi:= (id \otimes \Phi)\Delta_1$ and $\hat{\tilde{\delta}}_2 := \tilde{\Delta}_2$ over $C_2$, such that $\hat{\tilde{\delta}}_2 \Phi^{-1}:= (\Phi^{-1} \otimes id)\tilde{\Delta}_2$.
Therefore, the codipterous coalgebra $(\D, \Delta^*, \delta_1)$ and the anti-codipterous coalgebra $(\D, \tilde{\Delta}^*, \hat{\tilde{\delta}}_2)$ are entangled since 
$(id \otimes \hat{\tilde{\delta}}_2)\delta_1 = (\delta_1 \otimes id )\hat{\tilde{\delta}}_2$ is verified.
Let us define the maps
$\overrightarrow{d}: \ \D \xrightarrow{} \D \otimes \D$, such that $\overrightarrow{d}
=\Delta^*  - \delta_1$ and
$\overleftarrow{d}: \ \D \xrightarrow{} \D \otimes \D$, such that $\overleftarrow{d}
=    \tilde{\Delta}^*  - \hat{\tilde{\delta}}_2$.
We get:
$(id \otimes \overleftarrow{d})\overrightarrow{d}= (\overrightarrow{d} \otimes id)\overleftarrow{d}.$ Moreover, $(\overrightarrow{d} \otimes id) \Delta^* = (id \otimes \Delta^*)\overrightarrow{d}$ and $(id \otimes \overleftarrow{d} ) \tilde{\Delta}^* = (\tilde{\Delta}^* \otimes id)\overleftarrow{d}$.
\end{theo}
\Proof
Straightforward.
\eproof
\begin{defi}{[bimodule]}
Let the codipterous coalgebra $(\D, \Delta^*, \delta_1)$ and the anti-codipterous coalgebra $(\D, \tilde{\Delta}^*, \hat{\tilde{\delta}}_2)$ be entangled
as in the previous theorem \ref{dipt4}. Suppose all the involved coproducts
are unital homomorphisms. 
As bridges are coassociative, we embed $\D ^{\otimes 2}$
into a {\it{$\D$-bimodule}} by defining the following products: 
Let $c, x, y \in \D$,
\begin{eqnarray*}
x \ {\hat{\tilde{\circ}}} \ \overleftarrow{d}(c) &=& \hat{\tilde{\delta}}_2(x)\overleftarrow{d}(c); \ \ \ \ \overleftarrow{d}(c) \ {\hat{\tilde{\circ}}} \ y = \overleftarrow{d}(c) \hat{\tilde{\delta}}_2(y).
\end{eqnarray*}
\begin{eqnarray*}
x \circ \overrightarrow{d}(c) &=& \delta_1(x)\overrightarrow{d}(c); \ \ \ \ \overrightarrow{d}(c) \circ y = \overrightarrow{d}(c)\delta_1(y).
\end{eqnarray*}
\end{defi}
\begin{theo}
Let the codipterous coalgebra $(\D, \Delta^*, \delta_1)$ and the anti-codipterous coalgebra $(\D, \tilde{\Delta}^*, \hat{\tilde{\delta}}_2)$ be entangled
as in the previous theorem \ref{dipt4}. Suppose all the involved coproducts
are unital homomorphisms. 
Then $\overleftarrow{d},\overrightarrow{d}$ are Ito derivatives.
\end{theo}
\Proof
Let $x,y \in \D$.
We have $\overleftarrow{d}(1) =0= \overrightarrow{d}(1)$. Moreover,
$\overrightarrow{d}(x)\overrightarrow{d}(y) = \Delta_*(xy) + \delta_1(xy) - \Delta_*(x)\delta_1(y) - \delta_1(x)\Delta_*(y)$, i.e.
$ \overrightarrow{d}(xy)= \overrightarrow{d}(x)\overrightarrow{d}(y) + \overrightarrow{d}(x) \circ y + x \circ \overrightarrow{d}(y).$
Similarly,
$ \overleftarrow{d}(xy)= \overleftarrow{d}(x)\overleftarrow{d}(y) + \overleftarrow{d}(x) \hat{\tilde{\delta}}_2(y) + \hat{\tilde{\delta}}_2(x)\overleftarrow{d}(y)$, that is $ \overleftarrow{d}(xy)= \overleftarrow{d}(x)\overleftarrow{d}(y) + \overleftarrow{d}(x) \ {\hat{\tilde{\circ}}} \ y + x \ {\hat{\tilde{\circ}}} \ \overleftarrow{d}(y)$.
\eproof
\begin{exam}{$[Sl_q(2) \dpi \widetilde{Sl_q(2)}]$ } 
In \cite{codialg1}, the $(4,1)$-De Bruijn graph,
\begin{center}
\includegraphics*[width=3.5cm]{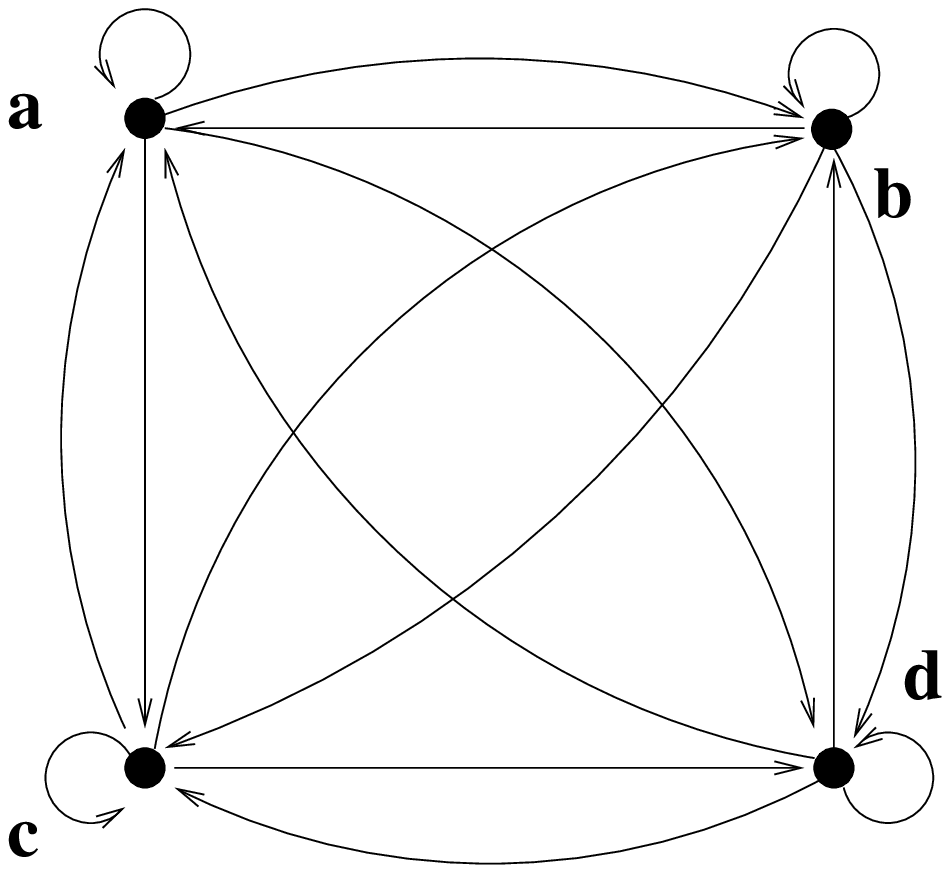}

\textbf{\begin{scriptsize} The $(4,1)$-De Bruijn graph. \end{scriptsize}}
\end{center}
was tiled with the geometric supports of two coassociative coalgebras, defining an achiral $L$-coalgebra, the coassociative coalgebra
represented by $\F$ and the coassociative coalgebra represented by the directed graph
defined by  $\tilde{\Delta}$ verifying $
\tilde{\Delta}b = b \otimes b + a \otimes d, \
\tilde{\Delta}c = c \otimes c + d \otimes a, 
\tilde{\Delta}a = b \otimes a + a \otimes c, \ 
\tilde{\Delta}d = c \otimes d + d \otimes b.$

Suppose $E$ is an algebra with unit 1, generated by eight elements $a,b,c,d,x,y,z$ and $u$.
Let $q$ be an invertible scalar of $k$. Recall that the Hopf algebra $Sl_q(2)$ is generated as an algebra by $a,b,c,d$ obeying the algebraic relations: $ba=qab, \ ca =qac, \ bc=cb, \ dc=qcd, \ db=qbd,
ad-da = (q^{-1}-q)bc, ad - q^{-1}bc =1$. The coproduct $\Delta_1$ verifies
$
\Delta_1 a = a \otimes a + b \otimes c, \ 
\Delta_1 b = a \otimes b + b \otimes d, \ 
\Delta_1 c = d \otimes c + c \otimes a, \ 
\Delta_1 d = d \otimes d + c \otimes b.
$
The antipode map $S$ is such that
$S_1(a) =d, \ S_1(d) =a, \ S_1(b) =-qb, \ S_1(c) = -q^{-1}c.$ Similarly,
Define $\widetilde{Sl_q(2)}$, the Hopf algebra generated by $x,y,z,u$, by the following coproduct:
$
\tilde{\Delta}_2 y = y \otimes y + x \otimes u, \ 
\tilde{\Delta}_2 z = z \otimes z + u \otimes x, \
\tilde{\Delta}_2 x = y \otimes x + x \otimes z, \ 
\tilde{\Delta}_2 u = z \otimes u + u \otimes y,$
with directed graph:
\begin{center}
\includegraphics*[width=3.5cm]{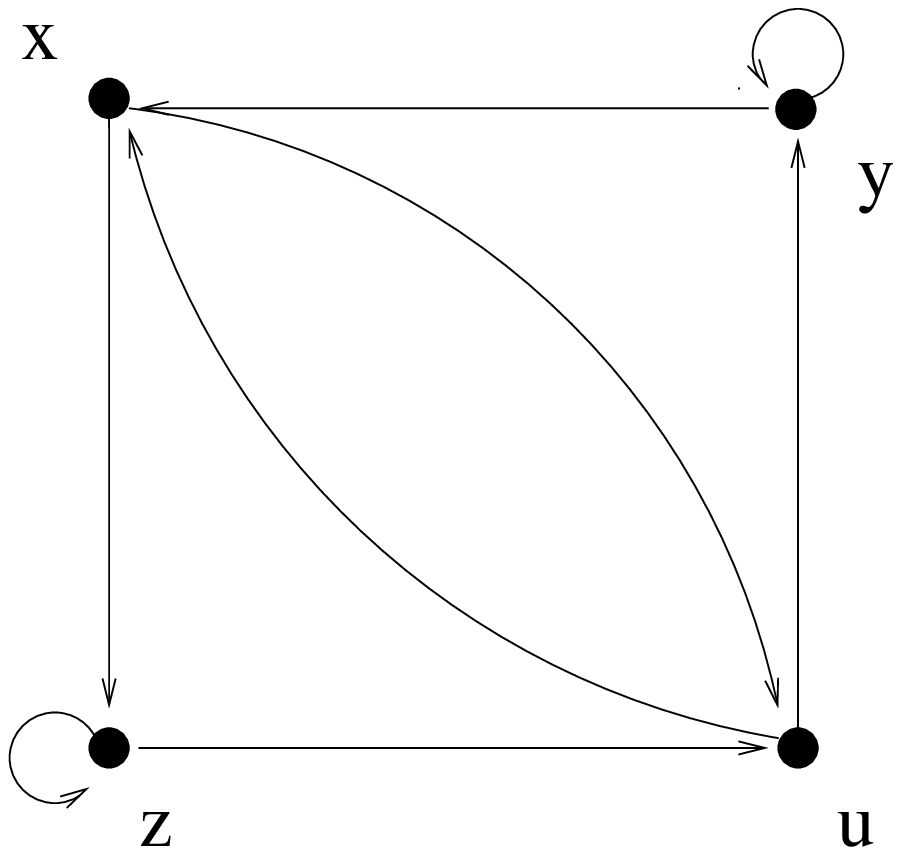}

\textbf{\begin{scriptsize} Geometric support of $\widetilde{Sl_q(2)}$. \end{scriptsize}}
\end{center}
The algebraic relations are $xy=qyx, \ uy = qyu, \ xu=ux, \ zu= quz, \ zx=qxz, \ yz-zy=(q^{-1}-q)xu, \ yz-q^{-1}xu =1$. The antipode map
is defined by $\tilde{S}_2(y) =z, \ \tilde{S}_2(z) =y, \ \tilde{S}_2(x) =-qx, \ 
\tilde{S}_2 (u) = -q^{-1}u.$
Define, as a channel map, the homomorphism $M: E \xrightarrow{} E$ such that
$M(a)=y, \ M(b)=x, \ M(c)=u, \ M(d)=z$. Notice that $\tilde{\Delta}_2 := (M \otimes M) \Delta M^{-1} $.
The bridges are defined by $\delta_1(x) := a \otimes x + b \otimes z, \ 
\delta_1(y) := a \otimes y + b \otimes u, \ \delta_1(z) := d \otimes z + c \otimes x, \ \delta_1(u) := d \otimes u + c \otimes y$, on $\widetilde{Sl_q(2)}$ and by $\Delta_1$ over $Sl_q(2)$. Similarly,
$\tilde{\delta}_2 (a) := y \otimes a + x \otimes c, \ 
\tilde{\delta}_2(b) := y \otimes b + x \otimes d, \ \tilde{\delta}_2(c) := z \otimes c + u \otimes a, \ \tilde{\delta}_2(d) := z \otimes d + u \otimes b$, over $Sl_q(2)$ and by $\tilde{\Delta}_2$ over $\widetilde{Sl_q(2)}$.
\begin{prop}
Define the coproduct $\Delta_* =\Delta_1$ over $Sl_q(2)$ and 
$\Delta_* =\tilde{\Delta}_2$ over $\widetilde{Sl_q(2)}$.
Equipped with the bridges $\delta_1$ for $Sl_q(2)$ and $\tilde{\delta}_2$, the codipterous coalgebra $(Sl_q(2), \Delta_*, \delta_1)$ is entangled with the codipterous coalgebra $(\widetilde{Sl_q(2)}, \Delta_*, \tilde{\delta}_2)$. Moreover
the bridges are unital homomorphisms.
\end{prop}
\Proof
The bridges are unital homomorphisms comes from the homomorphism property
of $M$.
\eproof

\noindent
The space $[Sl_q(2) \dpi \widetilde{Sl_q(2)}]$ has the following geometric support:
\begin{center}
\includegraphics*[width=7cm]{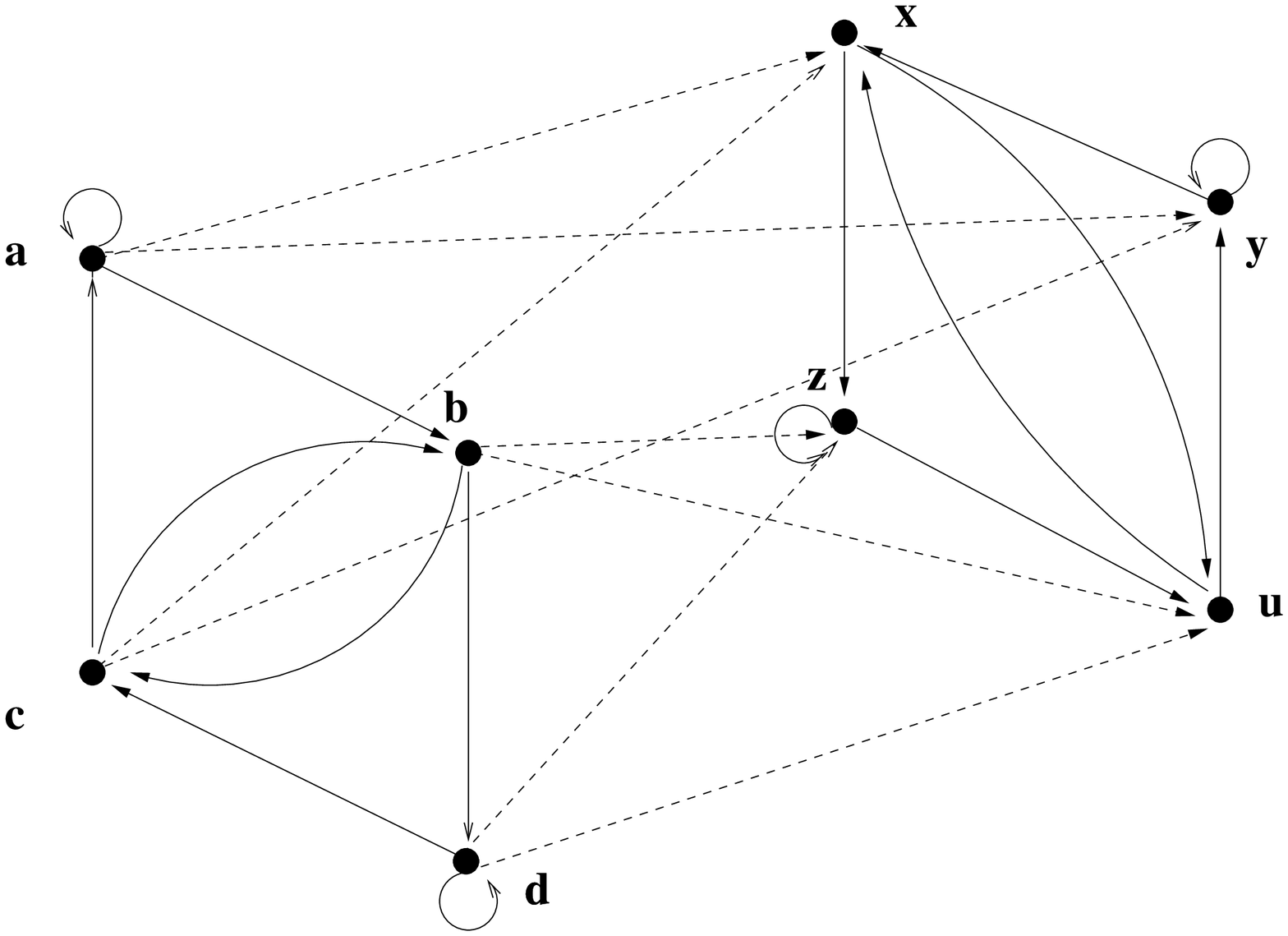}

\textbf{\begin{scriptsize} Geometric support of $[Sl_q(2) \dpi \widetilde{Sl_q(2)}]$. (Among the bridges, only $\delta_1$ is represented.) \end{scriptsize}}
\end{center}
\end{exam}
\begin{exam}{$[SU_2(q) \dpi \widetilde{SU_2(q)}]$}
The same results hold for $SU_2(q)$. Indeed, let $E$ be the algebra generated by $a, \ a^*, \ c, \ c^*, \ x, \ z^*, \ z, \ x^*$.
Recall that the Hopf algebra $SU_2(q)$, generated by $a, \ a^*, \ c, \ c^* $, whose coassociative coproduct $\Delta_1$ is described by
$\Delta_1 a = a \otimes a - qc^* \otimes c$ and $\Delta_1 c = c \otimes a + a^* \otimes c$ is in fact a part of an achiral $L$-coalgebra whose the other
coproduct \cite{codialg1} is described by  
$\tilde{\Delta}_1 a := c^* \otimes a + a \otimes c$ and
$\tilde{\Delta}_1 c := c \otimes c - q^{-1} a^* \otimes a$.
Similarly, define, as a channel, the $*$-homomorphism $M$ such that
$M(a) =z^*, \ M(c^*) = -q^{-1}x$.
Define also $\widetilde{SU_2(q)}$, generated by $\ x, \ z^*, \ z, \ x^*$, such that $\tilde{\Delta}_2 := (M \otimes M) \Delta_1 M^{-1}$
we can also construct the bridges $\delta_1$ and $\tilde{\delta_2}$, which are unital homomorphisms
and define the entanglement $[SU_2(q) \dpi \widetilde{SU_2(q)}]$.
\end{exam}
\Rk
In \cite{codialg1}, we have shown that the $(n^2,1)$-De Bruijn graph was tiled by (the geometric supports of) $n$ coassociative coalgebras. We can create $n$ transformations, via channel maps of $Sl_n(q)$ and glue them together. The space so obtained will have $n$ boundaries.
\section{Entanglement with a coassociative Markov $L$-coalgebra}
\subsection{Entanglement with a $(n,1)$-De Bruijn codialgebra}
In \cite{codialg1}, we have constructed $L$-cocommutative coassociative codialgebras via coassociative Markov $L$-coalgebras which are for instance the $(n,1)$-De Bruijn graphs
and unital associative algebras. By $(D_{(n,1)}, (x_i)_{1 \leq i \leq n})$, 
we mean the coassociative  codialgebra generated by $(x_i)_{1 \leq i \leq n}$ as a $k$-vector space and whose coproducts are $\Delta_M x_i := x_i \otimes \Sigma$, where $\Sigma :=\sum_{1 \leq j \leq n} \ x_j$ and $\tilde{\Delta}_M x_i := \Sigma \otimes x_i$, for all $1 \leq i \leq n$. We call such a codialgebra the $(n,1)$-De Bruijn codialgebra.
\begin{prop}
\label{dipt5}
Fix $n>0$.
Let $E := G \oplus C$ be a $k$-vector space, with $\dim E = 2n$, with basis
$((x_i)_{1 \leq i \leq n})$, $(a_i)_{1 \leq i \leq n})$, where
$(G := D_{(n,1)}, \ (x_i)_{1 \leq i \leq n}))$ is the $(n,1)$-De Bruijn codialgebra and $(C, \Delta)$ is a coassociative coalgebra, such that
$\Phi: G \xrightarrow{} C$ is a channel. Define the bridges $\delta_M$ (resp. $\tilde{\delta}_M$), such that $\delta_M := \Delta_M$, (resp. $\tilde{\delta}_M := \tilde{\Delta}_M$) over $G$ and $\delta_M \Phi:= (id \otimes \Phi)\Delta_M$, (resp. $\tilde{\delta}_M \Phi:= (id \otimes \Phi)\tilde{\Delta}_M$). Define the bridge $\delta$ verifying $\delta:= \Delta$ over $C$ and $\delta \Phi^{-1}:= (id \otimes \Phi^{-1})\Delta$. Define the
coproduct $\Delta_*$, by $\Delta_* = \Delta_M$ over $G$
and $\Delta_* = \Delta$ over $C$. We get
$(id \otimes \delta_M)\delta = (\delta \otimes id)\delta_M$,
i.e. the two codipterous coalgebras $(C \oplus G, \delta, \Delta_*)$ and $(C \oplus G, \delta_M, \Delta_*)$ are entangled. Such an entanglement is denoted by $[G \dpm C]$. Similarly, since $(\tilde{\delta}_M \otimes id)\delta = (id \otimes \delta)\tilde{\delta}_M$, by defining the
coproduct $\Delta^*$, such that $\Delta^* = \tilde{\Delta}_M$ over $G$
and $\Delta^* = \Delta$ over $C$, the two codipterous coalgebras $(C \oplus G, \delta, \Delta^*)$ and $(C \oplus G, \tilde{\delta}_M , \Delta^*)$ are entangled. Such an entanglement is denoted by $[\tilde{G} \dpmp C]$. 
\end{prop}
\Proof
Fix $n>0$ and denote by $G := D_{(n,1)}$ the $(n,1)$-De Bruijn codialgebra  generated by $(x_i)_{1 \leq i \leq n}$. Denote also by $(a_j)_{1 \leq j \leq n}$, the basis of $C$.
With the hypotheses and notation of the proposition \ref{dipt5}, we check that $(\Delta_* \otimes id)\delta_M =(\Delta_M \otimes id)\delta_M = (id \otimes \delta_M)\delta_M$ and that $(\Delta_* \otimes id)\delta =(\Delta \otimes id)\delta = (id \otimes \delta)\delta$. Let us prove the equality $(id \otimes \delta_M)\delta = (\delta \otimes id)\delta_M$. Recall that for 
all $i$, with $1 \leq i \leq n$, the Markovian coproduct $\Delta_M$ is
such that $\Delta_M x_i := x_i \otimes \Sigma$, where $\Sigma :=\sum_{1 \leq j \leq n} \ x_j$. Let $x_i \in G$, with $1 \leq i \leq n$, there exists an unique $c \in C$, such that $\Phi^{-1}(c) =x_i$. Using the Sweedler's notation, $\Delta c := \sum c_{(1)} \otimes c_{(2)}$, we get
$x_i \xrightarrow{\delta} \sum c_{(1)} \otimes \Phi^{-1}(c_{(2)}) \xrightarrow{(id \otimes \delta_M)} \sum c_{(1)} \otimes (\Phi^{-1}(c_{(2)}) \otimes \Sigma)$ and $x_i \xrightarrow{\delta_M} \sum x_i \otimes \Sigma \xrightarrow{(\delta \otimes id)} (\sum c_{(1)} \otimes \Phi^{-1}(c_{(2)})) \otimes \Sigma$. Similarly, we show that the same equality over $C$. The last assertion of the theorem is straightforward (recall that $\tilde{\Delta}_M := \tau \Delta_M$).
\eproof
\begin{exam}{$[D_{(4,1)} \dpm \F]$}
We have represented a part of the geometric support of $[D_{(4,1)} \dpm \F]$, where $D_{(4,1)}$ is the $(4,1)$-De Bruijn codialgebra generated as a $k$-vector space by $a',b',c',d'$.
\begin{center}
\includegraphics*[width=8cm]{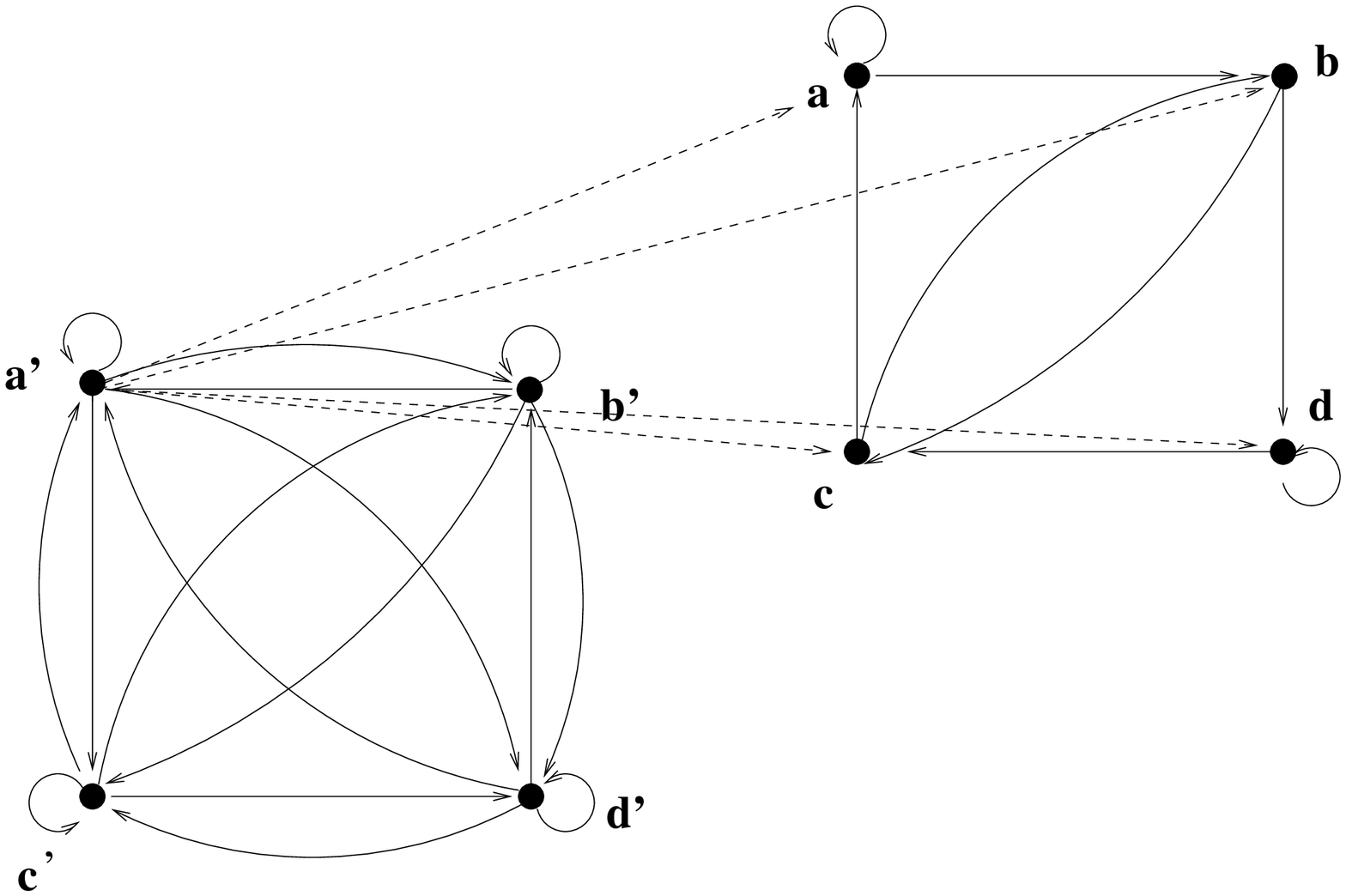}

\textbf{\begin{scriptsize} A part of the geometric support of $[D_{(4,1)} \dpm \F]$. \end{scriptsize}}
\end{center}
\end{exam}
\begin{exam}{$[D_{(1,1)} \dpm D_{(1,1)}]$}
Similarly, we can draw the self-entanglement of the $(1,1)$-De Bruijn codialgebra. Suppose the $k$-vector space $E$ is generated by two elements $X, Y$ and the channel map $\Phi': \ E \xrightarrow{} E$ such that $\Phi'(X) =Y$,
maps $(D_{(1,1)}, (X))$ into $(D_{(1,1)}, (Y))$.
\begin{center}
\includegraphics*[width=9cm]{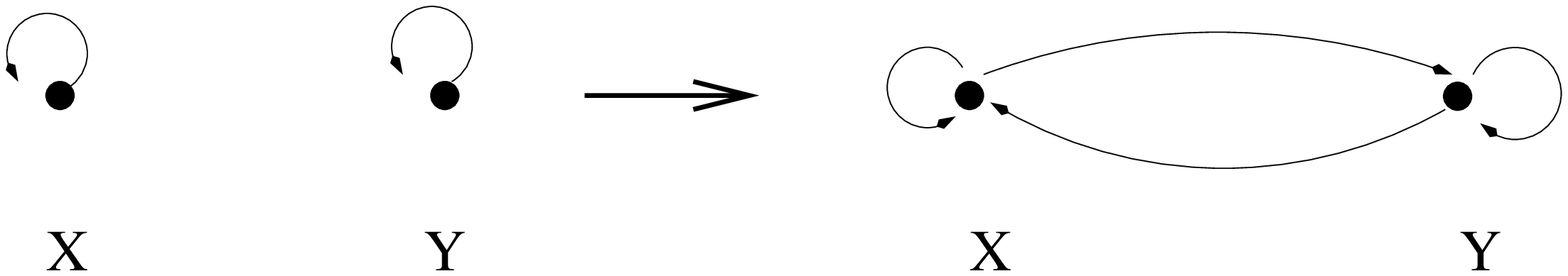}

\textbf{\begin{scriptsize} Right hand side: $\partial [D_{(1,1)} \dpm D_{(1,1)}]$. Left hand side: $[D_{(1,1)} \dpm D_{(1,1)}]$.\end{scriptsize}}
\end{center}
\end{exam}
Recall that $(D_{(1,1)}, (X))$ generated as a $k$-vector space, by $X$,
is defined by the coassociative coproduct $\Delta_M(X) := X \otimes X$. Therefore $(D_{(1,1)}, (Y))$, generated by $\Phi'(X):=Y$, as a $k$-vector space, is defined by the coassociative coproduct \footnote{By definition, a channel is a morphism of coalgebra.} $\Delta'_M(Y) := Y \otimes Y$, hence the geometric support
of the boundaries of $[D_{(1,1)} \dpm D_{(1,1)}]$. The bridge on  $(D_{(1,1)}, (X))$ is defined as $\delta_M:= \Delta_M$ over $(D_{(1,1)}, (X)$ and 
$\delta_M(Y) := \delta_M \Phi(X) := X \otimes \Phi'(X) := X \otimes Y$. Similarly, the bridge on  $(D_{(1,1)}, (Y))$ is defined as $\delta'_M:= \Delta'_M$ over $(D_{(1,1)}, (Y))$ and 
$\delta'_M(X) := \delta'_M \Phi^{' -1}(Y) :=  Y \otimes \Phi^{' -1}(Y) :=  Y \otimes  X$, hence the $(2,1)$-De Bruijn graph as geometric support of $[(D_{(1,1)}, \Delta_M, \delta_M) \dpm (D_{(1,1)}, \Delta'_M, \delta'_M)]$. This result is general. 
\begin{prop}
\label{bruj}
Fix $n >0$. Let $E := k\bra x_1, \ldots
x_n \ket \oplus k\bra y_1, \ldots
y_n \ket$ be a $k$-vector space, 
consider the $(n,1)$-De Bruijn codialgebra $(D_{(n,1)}, (x_j)_{1 \leq j \leq n})$. Consider a channel $\Phi': k\bra x_1, \ldots
x_n \ket \xrightarrow{} k\bra y_1, \ldots
y_n \ket$ such that $\Phi'(x_i)=y_i$ for all $i$. Then, the geometric support of $[D_{(n,1)} \dpm D_{(n,1)}]$ is the $(2n,1)$-De Bruijn graph.
\end{prop}
\Proof
Straightforward.
\eproof
\subsection{Entanglement with an unital associative algebra}
Let $A$ be an unital algebra with unit 1. We recall from \cite{Coa} that $A$ carries a non trivial Markov $L$-bi-dialgebra called the flower graph
with coproducts $\delta_f(a) = a \otimes 1$ and $\tilde{\delta}_f(a)= 1 \otimes a$, for all $a \in A$.
We call such a Markov $L$-coalgebra a {\it{flower graph}} because it is the concatenation of petals:
\begin{center}
\includegraphics*[width=4.5cm]{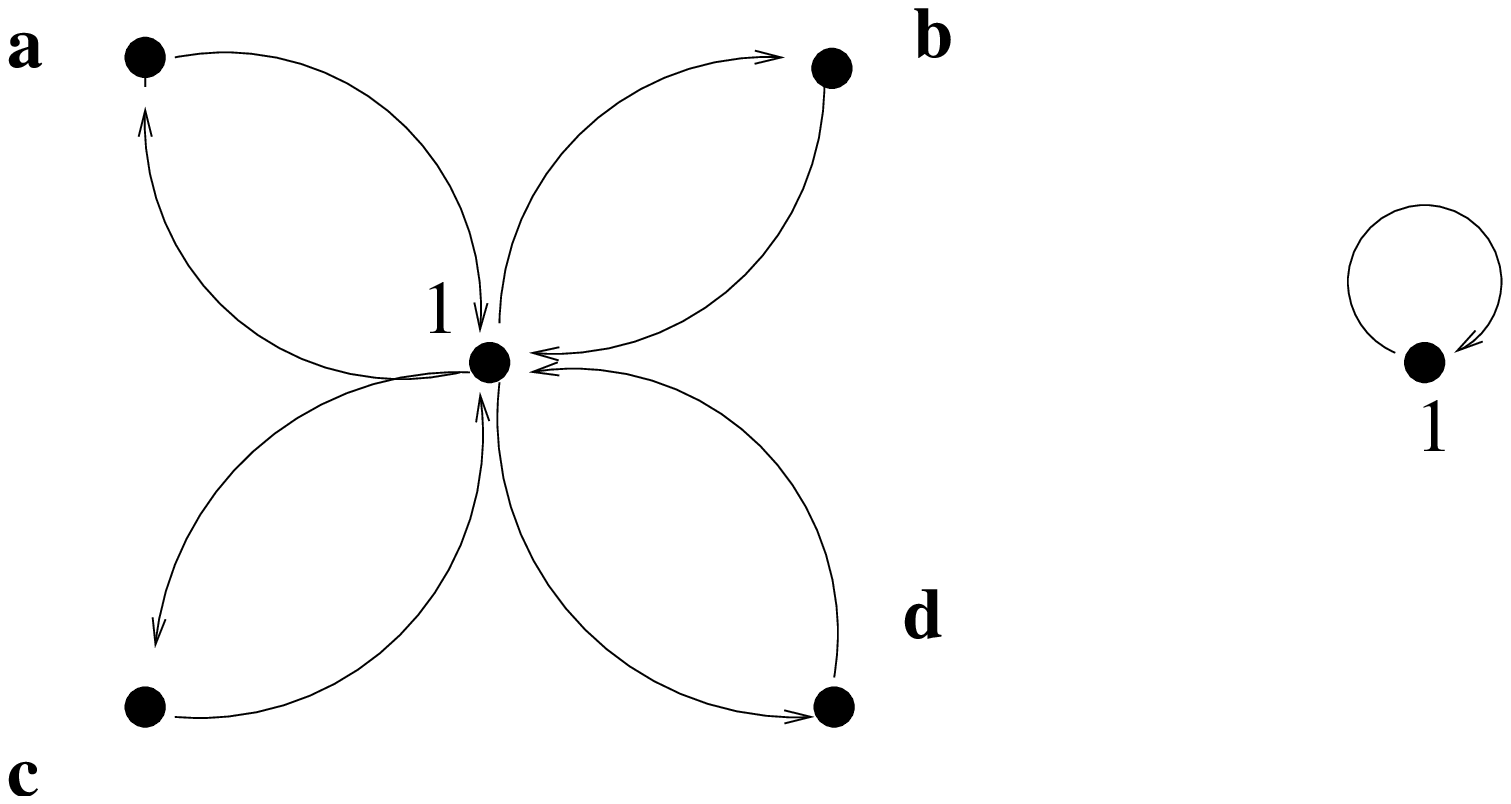}

\textbf{\begin{scriptsize} Geometric support of an algebra $A$. \end{scriptsize}}
\end{center}
\begin{prop}
Let $E$ be an algebra unit 1, composed by the bialgebra $(C, \Delta)$ and by an algebra
$(A,\tilde{\delta}_f, \delta_f)$, viewed as a $L$-bialgebra, (thus with the flower graph as geometric support). Suppose $A$ and $C$ are related by the homomorphism $\Phi: A \xrightarrow{} C$, as a chanel. Extend
the coproduct $\tilde{\delta}_f$ over the coassociative bialgebra $C$ by requiring
that $\tilde{\delta}_f \Phi = (id \otimes \Phi) \tilde{\delta}_f$. In addition, define the bridge $\delta$ such that $\delta := \Delta$ over $C$, and $\delta \Phi^{-1} = (id \otimes \Phi^{-1}) \Delta$. If $\Delta_* := \Delta_f:= \tilde{\delta}_f + \delta_f$ over $A$ and $\Delta_* := \Delta$ over $C$, then
the codipterous bialgebra $(A + C, \Delta_*, \tilde{\delta}_f)$ is entangled to the codipterous bialgebra $(A + C, \Delta_*, \delta)$, since
$(\tilde{\delta}_f \otimes id)\delta = (id \otimes \delta)\tilde{\delta}_f$. Such a chiral entanglement is denoted by $[\tilde{A} \dpmpf C]$. Similarly, since
$(id \otimes \delta_f)\delta = (\delta \otimes id )\delta_f$, $[A \dpmf C]$ is also a chiral entangled codipterous bialgebra.
\end{prop}
\Proof
Straightforward.
\eproof
\begin{exam}{$[ \tilde{A} \dpmpf Sl_q(2)]$} Let the algebra $E $ be generated by $a,b,c,d,a',b',c',d'$. Set $A:= k \bra a',b',c',d',1 \ket$. Suppose $k \bra a,b,c,d \ket$ stands for the Hopf algebra $Sl_q(2)$.
We represent here a face of the geometric support of $[\tilde{A} \dpmpf Sl_q(2)]$. The channel $\Phi: A \xrightarrow{} Sl_q(2)$ is such that $\Phi(a) := a', \Phi(b) := b', \Phi(c) := c', \Phi(d) := d'$.
\begin{center}
\includegraphics*[width=8cm]{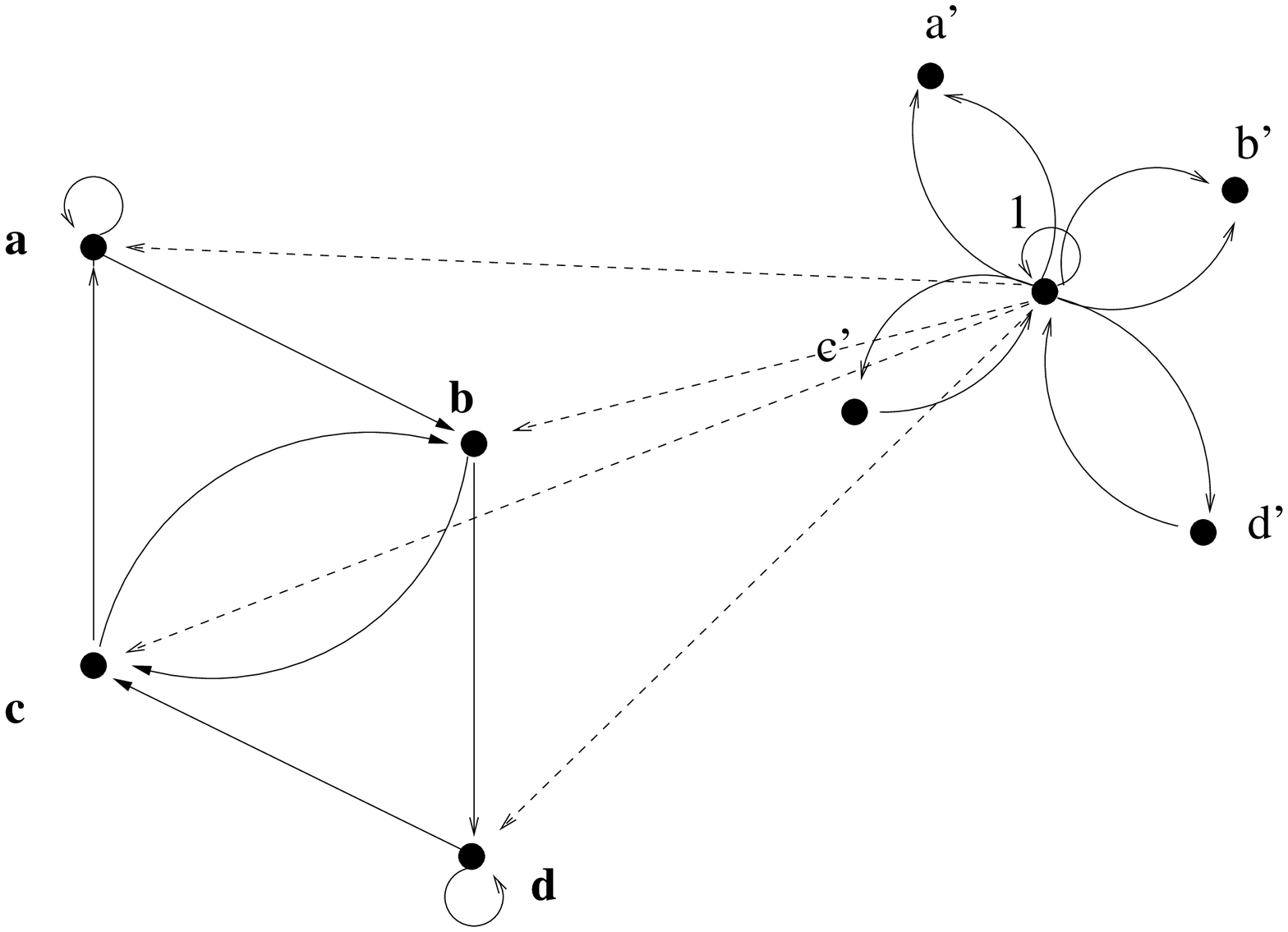}

\textbf{\begin{scriptsize} Face of the geometric support of $[ \tilde{A} \dpmpf Sl_q(2)]$. Among bridges, only $\tilde{\delta}_f$ is represented. \end{scriptsize}}
\end{center}
\end{exam}
We recall a result from \cite{Coa} inspired from a theorem from \cite{Hudson1}.
\begin{theo}
Let $(C, \Delta)$ be a bialgebra, then the coproducts
$\overrightarrow{d}:= \Delta - \delta_f$ and $\overleftarrow{d}:= \Delta - \tilde{\delta}_f$ are Ito derivatives and verify  $(id \otimes \overrightarrow{d}) \overleftarrow{d} =
(\overleftarrow{d} \otimes id)\overrightarrow{d}$.
\end{theo}
\subsection{Complex based on Ito coproducts of a bialgebra}
\textbf{Notation:} In this subsection, $\sum_i id \otimes \ldots \otimes \textrm{EXP} \otimes \ldots id$ means
that the expression EXP is placed at position $i$. Recall that $\Delta_f:=
\delta_f + \tilde{\delta}_f$.

Let $(C, \Delta)$ be a  bialgebra with unit 1. Let us establish a complex whose boundaries are based on the difference \footnote{Such a difference is used in the definition of connexe bialgebras.} $\Delta -\Delta_f$ and the Ito coproducts $\overrightarrow{d}:= \Delta - \delta_f$ and $\overleftarrow{d}:= \Delta - \tilde{\delta}_f$, where
$\Delta_f$ is the coassociative coproduct defined from the Markovian coproducts of the flower graph, naturally attached to $C$ as an algebra.
In this study, the primitive elements, i.e. $x \in C$ such that $\Delta(x) = x \otimes 1 + 1 \otimes x = \Delta_f(x)$ will vanish.
\begin{lemm}
Let us denote $ \partial'_0 = \Delta$ and for all $n>0$, $$\partial'_n =
\sum_{i=1}^{n+1} \ (-1)^{i+1} \ \underbrace{id \otimes \ldots \otimes id \otimes \Delta \otimes id \ldots \otimes id}_{n+1 \ \ \textrm{terms}}.$$
Fix $n>0$ and $a,a_1, \ldots, a_n \in C$. Then,
\begin{enumerate}
\item{$\partial'_{n+1} \tilde{\delta}_f (a_1) \otimes a_2 \otimes \ldots \otimes a_n =  \Delta(1) \otimes a_1 \otimes a_2 \otimes \ldots \otimes a_n - 1 \otimes \partial'_n (a_1 \otimes a_2 \otimes \ldots \otimes a_n).  $}
\item{$ (-1)^n \partial'_{n+1} a_1 \otimes a_2 \otimes \ldots \otimes \delta_f (a_n) = (-1)^n \partial'_n (a_1 \otimes a_2 \otimes \ldots \otimes a_n) \otimes 1 + a_1 \otimes a_2 \otimes \ldots \otimes a_n \otimes \Delta(1).$}
\item{$(\tilde{\delta}_f \otimes id)\tilde{\delta}_f = \Delta(1) \otimes id = (\Delta \otimes id) \tilde{\delta}_f $.}
\item{$ (id \otimes \delta_f)\delta_f = id \otimes \Delta(1) = (id \otimes \Delta)\delta_f$. }
\item{$ (\Delta \otimes id)\delta_f(a)= \Delta(a) \otimes 1, \ \ \ (id \otimes \Delta)\tilde{\delta}_f(a) = 1 \otimes \Delta(a) $.}
\end{enumerate}
\end{lemm}
\Proof
The proof is complete by noticing that $\Delta(1) = 1 \otimes 1$ and by using the definitions of
the coproducts $\delta_f(a) = a \otimes 1, \ \ \tilde{\delta}_f(a) = 1 \otimes a$, for all $a \in C$.
\eproof
\begin{theo}
\label{th}
Recall that $\overleftarrow{d} := \Delta - \tilde{\delta}_f$ and $\overrightarrow{d} := \Delta - \delta_f$. The sequence,
$$ 0 \xrightarrow{} C \xrightarrow{\Delta -\Delta_f} C^{\otimes 2} \xrightarrow{\partial_1 = \overleftarrow{d} \otimes id - id \otimes \overrightarrow{d}} C^{\otimes 3} \xrightarrow{\partial_2} C^{\otimes 4} \xrightarrow{\partial_3} \ldots$$
with $\forall n>0$,
$$\partial_n := \underbrace{\overleftarrow{d} \otimes id \otimes \ldots \otimes id}_{n+1 \ \ \textrm{terms}} + \sum_{j=2}^{n-1} \ \ (-1)^{j+1} \ id \otimes \ldots \otimes id \otimes
\Delta \otimes id \otimes \ldots \otimes id
+ (-1)^{n + 1}id \otimes id \otimes \ldots \otimes id \otimes \overrightarrow{d},$$ defines a complex.
The boundary operators verify:
\begin{enumerate}
\item {$\forall n>0 \ \ \partial_{n+1}\partial_n = 0$ and $\partial_1(\Delta -\Delta_f)=0$.}
\item{$\forall n>0 \ \ \partial_n (x_1, \ldots , x_n)$ is a multilinear map which is an Ito derivative in the first and last variables and a homomorphism in others variables.}
\end{enumerate}
\end{theo}
\Proof
We only have to prove the first item, since the second one comes from the very definition of the boundary operators.
From coassociativity coalgebra theory, we know that $ \partial'_{k+1} \partial'_k = 0, \ \ \forall k \in \mathbb{N}$.
Fix $n>0$.
\begin{eqnarray*}
\partial_{n+1} \partial_n &=& \partial'_{n+1} - \tilde{\delta}_f \otimes id \ldots \otimes id +(-1)^{n+1} id \otimes \ldots \otimes id \otimes \delta_f)
(\partial'_n \\
& & - \tilde{\delta}_f \otimes id \ldots \otimes id + (-1)^{n} id \otimes \ldots \otimes id \otimes \delta_f) \\
&=& \partial'_{n+1}\partial'_n - \partial'_{n+1}(\tilde{\delta}_f \otimes id \ldots \otimes id) + (-1)^{n} \partial'_{n+1}( id \otimes \ldots \otimes id \otimes \delta_f )\\
& &  - 1 \otimes \partial'_n + (\tilde{\delta}_f \otimes 1)\tilde{\delta}_f \otimes id \ldots \otimes id\\
& & +(-1)^{n+1} \partial'_n \otimes 1 - id \otimes \ldots \otimes id \otimes (id \otimes \delta_f)\delta_f.
\end{eqnarray*}
The equality $\partial_1(\Delta -\Delta_f) = 0$ is straightforward.
\eproof

\noindent
We can also express the boundary operators only in terms of $\Delta -\Delta_f$.
\begin{theo}
For all $n > 0$, $\partial_n = \sum_{i=1}^{n + 1} \ (-1)^{n+1} \ id \otimes \ldots \otimes (\Delta - \Delta_f) \otimes  \ldots \otimes id$.
\end{theo}
\Proof
The proof is straightforward by noticing that $id \otimes \tilde{\delta}_f = \delta_f \otimes id$.
\eproof
\begin{exam}{}
Consider the geometric support of $Sl_2(q)$. The difference $\Delta -\Delta_f$, at $a$ yields:
\begin{center}
\includegraphics*[width=10cm]{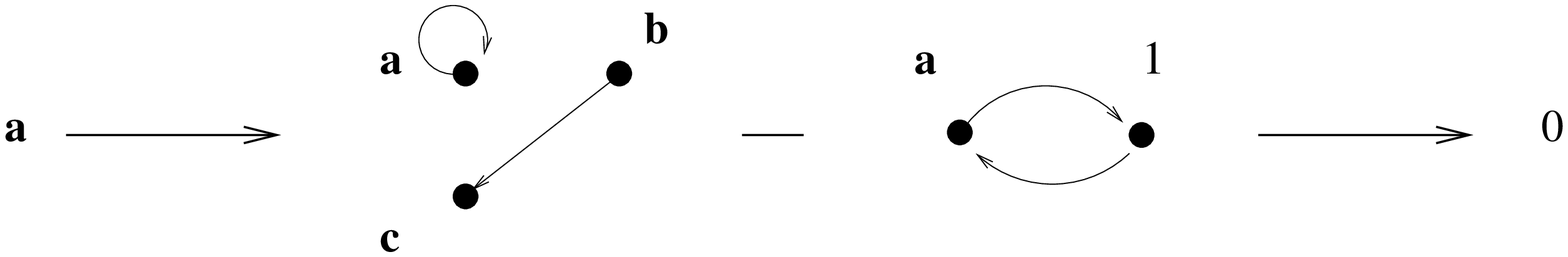}
\end{center}
and vanishes when the operator $\partial_1 = \overleftarrow{d} \otimes id - id \otimes \overrightarrow{d}$ is applied.
\end{exam}
\Rk
Notice that all the results of this subsection hold
for a coassociative coalgebra equipped with a group-like element $e$, by removing $1$ by $e$. Of
course, the coproducts $\overleftarrow{d}, \overrightarrow{d}$ are no longer Ito derivatives.
\section{Towards coassociative manifolds and $L$-molecules}
Thanks to the notion of entanglement of codipterous coalgebras and anti-codipterous coalgebras, we construct more and more
complicated directed graphs as well as a generalisation of coassociative coalgebras, bialgebras and Hopf algebras
by considering them
as simple structures, reminiscent of atoms in physics. The entanglement of atoms leads to molecules. We define:
\begin{defi}{[$L$-molecule]} 
Let $E$ be a $k$-vector space.
A {\it{$L$-molecule}} in $E$ is an entanglement of several codipterous coalgebras and / or anti-codipterous coalgebras constructed from coassociative coalgebras via channel maps.
\end{defi}
\begin{exam}{}
For instance, $[Sl_q(2) \dpe \widetilde{Sl_q(2)} \dpe Sl_q(2)]$ or
$[Sl_q(2) \dpe Sl_q(2) \dpe \cdots \dpe Sl_q(2)]$.
\end{exam}
By dualizing all these definitions, we can construct algebras with several laws and can envisage to study their associated operads.

\noindent
\textbf{Conjecture:}
Any operad, associated with a dualisation of a $L$-molecule is Koszul.

Let us end with a geometric interpretation of the notion of coproduct of a $L$-molecule $M$, generated as a $k$-vector space by a basis $(v_i)_{1 \leq i \leq \dim M}$. Recall that the family $(v_i)_{1 \leq i \leq \dim M}$ plays the r\^ole of the vertex set and $v_i \otimes v_j$ is symbolized
by a directed arrow $v_i \xrightarrow{} v_j$. In addition,
by the intersection of two directed graphs, we mean the intersection of their arrow sets.

\begin{defi}{[Coassociative covering]}
Let $W$ be a $k$-vector space, with basis $(v_i)_{1 \leq i \leq \dim W}$. 
Consider a  subspace $\mathcal{G}$ of $W^{\otimes 2}$.
Denote by $G_\delta$, the directed graph defined by the coassociative coproduct $\delta$. Let $\delta_1, \ldots, \delta_n: \ W \xrightarrow{} \mathcal{G}$ be $n$ coassociative coproducts. The directed graph $G$, associated with $\mathcal{G}$ has a 
{\it{coassociative covering}}  if
$\cup G_{\delta_{i}} = G$ and $\delta_i = \delta_j$ over $G_{\delta_{i}} \cap G_{\delta_{j}}$. Such a space is called a {\it{coassociative manifold}}.
\end{defi}
\begin{theo}
Any geometric support of a $L$-molecule is a coassociative manifold.
\end{theo}
\Proof
Straightforward by the very definition of a $L$-molecule.
\eproof

\noindent
Let us show now how to embed any non directed graph
into a coassociative manifold.
Let $G=(G_0, G_1)$ be a non directed graph, such that given two vertices $u,v$ of $G_0$, they are linked either with no edge or with an unique edge of $G_1$. Let $G^\natural$ be the directed graph such that
the vertex set $G^\natural=G_0$. The arrow set
$G^\natural _1$ is defined as follow: each
edge, $\bullet --- \bullet$, of $G_1$ which is not a loop, is removed by a bi-directed arrow $\bullet \leftrightarrows \bullet$. In addition, put a loop on each vertex without loop in $G$. 
\begin{center}
\includegraphics*[width=3cm]{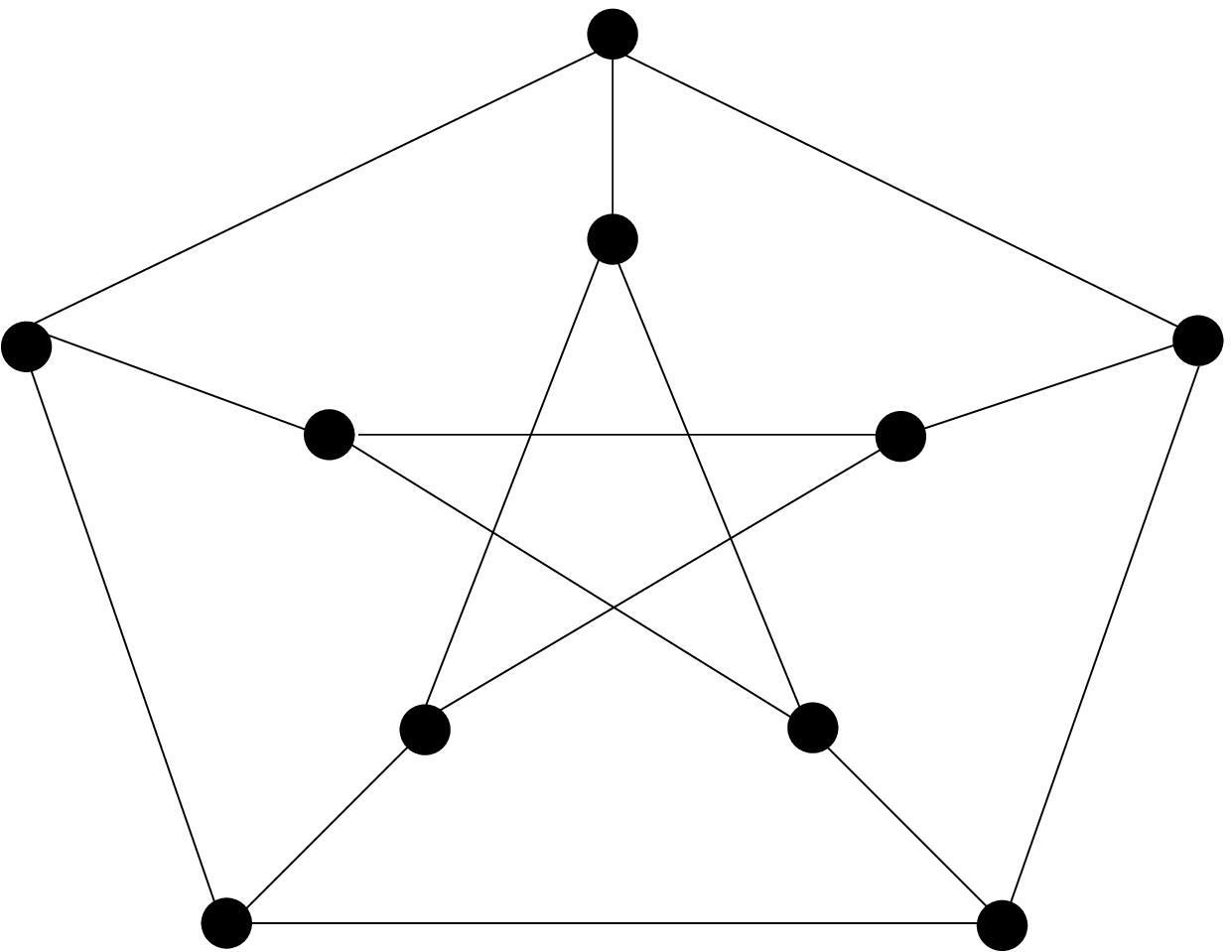}

\textbf{\begin{scriptsize} Example of non directed graph: The Petersen graph. \end{scriptsize}}
\end{center}
\begin{theo}
Let $G=(G_0, G_1)$ be a non directed graph such that two vertices $u,v$ of $G_0$, are linked either with no edge or with an unique edge of $G_1$. Then
$G^\natural$ is the geometric support of a $L$-molecule and thus of a coassociative manifold. \end{theo}
\Proof
Embed $G^\natural_0$ into its free $k$-vector space $kG^\natural_0$. Define the coassociative coproduct $\Delta_l v:= v \otimes v$ for all $v \in G^\natural_0$. Let $v,w \in G^\natural_0$, with $v \not= w$. If
$v \longrightarrow w$ is an arrow of $G^\natural_1$,
define the chanel map $\phi_{vw}: k\bra v \ket \xrightarrow{} k\bra w \ket$ such that $\phi_{vw}(v)=w$. Define the
bridges $\delta_{vw}w:= v \otimes w$ and $\delta_{vw}v:= v \otimes v = \Delta_l v$. Then
the family of bridges $(\delta_{s(a)t(a)})_{a \in G^\natural_1}$ is a family of coassociative coproducts. 
Moreover, $(kG^\natural_0, \Delta_l, \delta_{s(a)t(a)})_{a \in G^\natural_1}$ is a family of codipterous coalgebras and for all $v,w \in G_1$, 
with $v \not= w$, the space $(k\bra v \ket \oplus k\bra w \ket, \Delta_l, \delta_{vw}, \delta_{wv})$ is (chiral)
entangled, according to section 2. Therefore, the directed graph $G^\natural$ is covered by geometric supports 
of (Markovian) coassociative coalgebras. As the common intersection between two geometric supports is either
a loop or the empty set, $G^\natural$ is then a coassociative manifold since the coproducts are equal to $\Delta_l$ on a loop. 
\eproof
\NB
Notice also that, ``locally'' two vertices related
by a bi-directed arrow looks like a $(2,1)$-De Bruijn graph.
Moreover,
the map $\epsilon: kG^\natural_0 \xrightarrow{} k$
such that $\epsilon(v) =1$ for all $v \in G^\natural_0$ is a left counit.
\section{Conclusion}
In this ``coassociative geometry'', the r\^ole of all the coproducts of a coassociative manifold can be better understood. As a comparison, a topological manifold can be covered by
open sets, each one being homeomorphic to an open set of $\mathbb{R}^n$. Here, the r\^ole of an open set is played by a coassociative coalgebra, being a side of the $L$-molecule $M$. Denote by
$L(M, A)$, the space of linear mappings defined over $M$, with values in an algebra $A$, on such side or ``open set'' of $M$, there will be an
unique convolution product. By ``travelling'' over such a manifold, we will change
progressively the convolution product since $\delta_i=\delta_j$ over $G_{\delta_{i}} \cap G_{\delta_{j}}$. As an example, recall that the self-entanglement of $\F$, denoted by $M$, yields two pre-dendriform coalgebras, each one having three coproducts. Consider $L(M, A)$, where $A$ is an associative algebra.
Suppose we start with the side $\delta_1$ of $M$. This means that we live in an associative algebra $(L(M, A), \vdash)$ whose convolution product of two linear maps $f,g$, $f \vdash g$ is defined by $\delta_1$. Continuing our travel, suppose we arrive on the side $(\F, \Delta_1)$. As $\delta_1 := \Delta_1$ on the side $(\F, \Delta_1)$, we will obtain $f \vdash g := f \perp g \ (= f \dashv g )$. Leaving the side $(\F, \Delta_1)$ to go to the side
$\hat{\delta}_1$, we will arrive in an associative algebra $(L(M, A), \dashv)$ whose convolution product is now determined by $\hat{\delta}_1$.
Therefore, turning around the sides of the associative trialgebra $(L(M, A), \vdash,\dashv, \perp)$ leads us to pass ''progressively`` from $(L(M, A), \vdash)$ to $(L(M, A), \perp)$ to $(L(M, A), \dashv)$. Besides this dynamical point of view,
the formalism developed so far could have several applications in graph theory but also in non commutative stochastic processes. Quantum Levy processes on bialgebras are now well understood, see \cite{Schurmann,Franz} for instance. Would it be possible to develop such a theory on coassociative manifolds, or quantum Levy processes on bialgebras with several parameters, and thus with several coproducts? This work will pursued in \cite{Lerdiff}.

\noindent
\textbf{Acknowledgments:}
The author wishes to thank D. Petritis for fruitful advice for the redaction of this paper as well as J-L. Loday and U. Franz and M. Schurmann for useful discussions.

\bibliographystyle{plain}
\bibliography{These}

\end{document}